\documentclass[preprint,review,12pt]{elsarticle}

\usepackage{amsbsy,bm,amsmath,amssymb}
\usepackage{graphicx}
\usepackage{epsfig}

\newcommand{\nn}{\nonumber}
\newcommand{\bpsi}{\boldsymbol{\psi}}
\newcommand{\bphi}{\boldsymbol{\varphi}}

\newcommand{\mx}{\mathbf}

\newcommand{\rmax}{r_\mathrm{max}}

\journal{Computer Physics Communications}

\begin{document}

\begin{frontmatter}

\title{Orthogonal fast spherical Bessel transform on uniform grid}

\author{Vladislav V. Serov}

\cortext[author] {Corresponding author.\\\textit{E-mail address:} vladislav\_serov@mail.ru}
\address{
Department of Theoretical Physics, Saratov State University, 83
Astrakhanskaya, Saratov 410012, Russia}

\begin{abstract}
We propose an algorithm for the orthogonal fast discrete spherical Bessel transform on an uniform grid. Our approach is based upon the spherical Bessel transform factorization into the two subsequent orthogonal transforms, namely the fast Fourier transform and the orthogonal transform founded on the derivatives of the discrete Legendre orthogonal polynomials.
The method utility
is illustrated by its implementation for the numerical solution of the three-dimensional time-dependent Schr\"odinger equation.
\end{abstract}



\end{frontmatter}

\section{Introduction}

The discrete spherical Bessel transform (DSBT) arises in a number of applications, such as, e.g., the analysis of the cosmic microwave background \cite{Hamilton2000}, the numerical solution of the differential equations \cite{Bisseling1985,Kouri1987,Ronen2006}, and the numerical evaluation of multi-center integrals \cite{Talman2003,Toyoda2009}. Many different SBT algorithms have been proposed so far \cite{Talman1978,Sharafeddin1992,Lemoine1994,Toyoda2010}. But none of them possess all of the advantages of their trigonometric progenitor, namely the fast Fourier transform (FFT). These advantages are the performance fastness, the uniform coordinate grid, and the orthogonality.

An example of the problem requiring the simultaneous presence of all the advantages is the solving of the Schr\"odinger-type equation (SE) by means of the pseudospectral approach \cite{Lemoine1994}.
The grid uniformity provides the same accuracy of the wave function description in the whole domain of definition.
The grid identity for all orders of a spherical Bessel functions (SBF) allows to switch to the discrete variable representation (DVR) \cite{Lemoine1994}. 
The DSBT orthogonality is needed to provide the hermiticity of the radial part of the Laplacian operator in DVR. The lack of the Laplacian operator hermiticity impedes the convergence of iterative methods (such as conjugate gradient method) for the solution of matrix equations (which are obtained by DVR from the stationary SE). In the case of time-dependent SE, the hermiticity of the Laplacian operator is crucial for the conservation of the wave function norm during the time evolution. 

A pioneering approach based upon the convolution integral \cite{Talman1978,Siegman1977,Talman2009} requires a number of operations of the order of $N\log_2N$ for its performing, just like the FFT does, that means that it is quite fast. However it employs a strongly nonuniform grid (a node location exponentially depending on its number). Hence the attempts of its utilization for the SE solving \cite{Bisseling1985,Kouri1987} ended in problems with the strong near-center localization of a wave function. A method rest on the spherical Bessel functions expansion over the trigonometric functions \cite{Sharafeddin1992} also appears to be quite fast (requiring as few as $(\ell+1)N\log_2N$ operations) and employs a uniform grid. But it is not orthogonal and has stability difficulties because of the singular factors in the spherical Bessel functions expansion over the trigonometric functions. Next, a Gauss-Bessel quadrature based technique suggested in \cite{Lemoine1994,Lemoine2003} is orthogonal and converges exponentially, but it is not fast (as the number of operations required scales as $N^2$) and needs an $\ell$-dependent grid. Nevertheless its fast convergence and the near-uniform grid motivated to apply it for a time-dependent Gross-Pitaevsky equation \cite{Ronen2006}. Finally, an approach rest on the SBF integral representation via Legendre polynomials, proposed in \cite{Toyoda2010}, appears to be fast, makes use of the uniform grid, but it is not orthogonal. 

In the present work we are proposing the algorithm for the DSBT that is orthogonal, fast, and it implies the uniform grid. Our approach is based upon the SBT factorization into the two subsequent transforms, namely the FFT and the discrete orthogonal Legendre polynomials derivatives based transform. 

The paper is organized as follows. In Section \ref{Sec:DSBT}, we develop the orthogonal fast DSBT on a uniform grid. Next, in Section \ref{Sec:Testing}, the proposed method is tested via the evaluation of the Gaussian atomic functions transform and also the DSBT basis functions comparison to the exact SBFs. In Section \ref{Sec:Example} the DSBT- and DVR-based approach (DSBT-DVR) for the time-dependent Schr\"odinger equation (TDSE) solving is suggested and examined. The approach efficiency is illustrated by treating of the problem of the Hydrogen molecular ion ionization by laser pulse. Finally, in Section \ref{Sec:Conclusion} we briefly discuss the obtained results as well as the prospects of DSBT and DSBT-DVR application.

\section{Development of the method}\label{Sec:DSBT}

\subsection{Basic formulation}

A typical problem involving the spherical Bessel transform (SBT) is the plane wave expansion of a three-dimensional function $\Psi(r,\theta,\phi)$. The expansion over the spherical harmonics yields a radius-dependent function 
\begin{eqnarray}
\Psi_{\ell m}(r) = \oint Y_{\ell m}(\theta,\phi) \Psi(r,\theta,\phi) d\Omega. \label{YProj}
\end{eqnarray}
If the function $\Psi(r,\theta,\phi)$ has no singularities, then $\Psi_{\ell m}(r\to 0)\sim r^{\ell}$.

Let us introduce the SBT as 
\begin{eqnarray}
c_\ell(k) = \sqrt{\frac{2}{\pi}} \int_0^\infty \chi_{\ell}(kr) \psi_\ell(r) dr \label{BesselProj}
\end{eqnarray}
Here we perform the function substitution $\psi_\ell(r)=r\Psi_{\ell m}(r)$ (a magnetic quantum number is not used further, therefore from this point on we omit it from the denotation for the sake of simplification), then execute the expansion over the functions 
\begin{eqnarray}
\chi_{\ell}(x) = x j_\ell(x), \label{krBessel}
\end{eqnarray} 
where $j_\ell(x)$ is a spherical Bessel function (SBF) of the first kind. The functions $\chi_{\ell}(kr)$ satisfy the normalization condition $\int_0^\infty \chi_{\ell}(kr)\chi_{\ell}(k'r) dr = (\pi/2)\delta_{kk'}$. The pre-integral factor in \eqref{BesselProj} is introduced in order to make the transform \eqref{BesselProj} unitary. 

The beginning of our derivation coincides with the one in the work \cite{Toyoda2010}. But unlike its authors we are going to aim at the factorization of the SBT into the two separate transformations, namely the FFT and also the additional orthogonal transform which we denote Fourie-to-Bessel transform (FtB).
The SBF may be presented as
\begin{eqnarray}
j_\ell(z) = \frac{1}{2i^\ell} \int_{-1}^{1} P_\ell(\eta) \exp(iz\eta) d\eta \label{BesselExp}
\end{eqnarray}
where $P_\ell(\eta)$ is the Legendre polynomial of $\ell$-th order. Upon substituting the latter expression into Eq.\eqref{BesselProj}, we obtain
\begin{eqnarray*}
c_\ell(k) = \sqrt{\frac{2}{\pi}} \frac{1}{2i^\ell} \int_{-1}^{1} P_\ell(\eta) \int_0^\infty kr e^{ikr\eta}\psi_\ell(r) dr d\eta
\end{eqnarray*}
Here the integral over $r$ is different from the Fourier transform of the function $\psi_\ell(r)$ by the presence of the integrand factor $kr$. This factor might be represented as a result of taking a derivative of $e^{ikr\eta}$ over $\eta$. Thus we get the expression
\begin{eqnarray*}
c_\ell(k) = \sqrt{\frac{2}{\pi}} \frac{1}{2i^{\ell+1}} \int_{-1}^{1} P_\ell(\eta) \frac{\partial}{\partial\eta}\int_0^\infty e^{ikr\eta}\psi_\ell(r) dr d\eta
\end{eqnarray*}
Making use of the Legendre polynomials parity condition $P_\ell(-\eta)=(-1)^\ell P_\ell(\eta)$, one may further reduce
 the integral over $\eta$ from $-1$ to $1$ to the one in the limits from $0$ to $1$ as
\begin{eqnarray}
c_\ell(k) = \sqrt{\frac{2}{\pi}} \frac{1}{2i^{\ell+1}} \int_{0}^{1} P_\ell(\eta) \frac{\partial}{\partial\eta}
\int_0^\infty [e^{ikr\eta}-(-1)^\ell e^{-ikr\eta}]\psi_\ell(r) dr
d\eta \label{c_vs_psi}
\end{eqnarray}

Next, let us define a new function
\begin{eqnarray}
\tilde{c}_\ell(k) = \sqrt{\frac{2}{\pi}} \frac{1}{2i^{\ell+1}}
\int_0^\infty [e^{ikr}-(-1)^\ell e^{-ikr}]\psi_\ell(r) dr \label{tilde_c_def}
\end{eqnarray}
The term $[e^{ikr}-(-1)^\ell e^{-ikr}]/(2i^{\ell+1})$ is equal to $(-1)^{\left\lceil\ell/2\right\rceil}\sin(kr)$ for the even $\ell$ and to $(-1)^{\left\lceil\ell/2\right\rceil}\cos(kr)$ for the odd ones. Hence the expression \eqref{tilde_c_def} appears to be correspondingly the sine/cosine Fourier transform of the function $\psi_\ell(r)$, depending on $\ell$ being even or odd.

In terms of the new denotation the formula \eqref{c_vs_psi} can be rewritten as
\begin{eqnarray*}
c_\ell(k) =  \int_{0}^{1} P_\ell(\eta) \frac{\partial\tilde{c}_\ell(k\eta)}{\partial\eta}  d\eta
\end{eqnarray*}
Upon making the substitution $\eta=q/k$ this expression takes the following form
\begin{eqnarray*}
c_\ell(k) &=& \int_{0}^{k} P_\ell(q/k) \frac{d\tilde{c}_\ell(q)}{dq} dq
\end{eqnarray*}
Let us perform the integration by parts, then move the derivative over $q$ to the Legendre polynomial. As a result we obtain the following formula for the FtB
\begin{eqnarray}
c_\ell(k) &=& \tilde{c}_\ell(k) - \int_{0}^{k} \frac{P_\ell'(q/k)}{k} \tilde{c}_\ell(q) dq \label{c_tildec} 
\end{eqnarray}
It is easily seen that $c_0(k)=\tilde{c}_0(k)$, just as expected, since $\chi_{0}(kr)=\sin(kr)$ and the Bessel expansion coincides with the Fourier expansion at $\ell=0$.

One may rewrite \eqref{c_tildec} in the operator form as $c_\ell(k)=\hat{T}\tilde{c}_\ell(k)$, where the integral transform operator $\hat{T}$ has the kernel
\begin{eqnarray}
T(k,q) &=& \delta(q-k) - \theta(k-q) \frac{P_\ell'(q/k)}{k} .
\end{eqnarray}
Here 
\begin{eqnarray}
\theta(x) = \left\{
\begin{array}{ll}
0, & x< 0; \\
1/2, & x=0; \\
1, & x > 0.
\end{array}\right. \label{step_func}
\end{eqnarray}
is the Heaviside step function.
The FtB operator $\hat{T}$ must be unitary (that is $\hat{T}^{-1}=\hat{T}^T$), hence the inverse transform is $\tilde{c}_\ell(k)=\hat{T}^Tc_\ell(k)$. The inverse FtB might be explicitly defined as
\begin{eqnarray}
\tilde{c}_\ell(k) &=& c_\ell(k) - \int_{k}^{\infty} \frac{P_\ell'(k/q)}{q} c_\ell(q) dq. \label{tildec_c} 
\end{eqnarray}

The substitution of \eqref{tildec_c} into \eqref{c_tildec} demonstrates that \eqref{tildec_c} is indeed inverse in respect to \eqref{c_tildec}, due to the condition
\begin{eqnarray}
&& \int_0^{\min(k_2,k_1)} \frac{P'_\ell(q/k_2)}{k_2} \frac{P'_\ell(q/k_1)}{k_1} dq = \nonumber \\ && \frac{P'_\ell(k_2/k_1)}{k_1}\theta(k_1-k_2) + \frac{P'_\ell(k_1/k_2)}{k_2}\theta(k_2-k_1). \label{intPP_cond}
\end{eqnarray}
This condition holds true since for any polynomial $p(x)$ of the order $s\leq\ell$ true is the expression
\begin{eqnarray}
&& \int_{-1}^1 P_\ell'(x) p(x) dx = p(1)-(-1)^{\ell}p(-1) \label{intPp_cond}
\end{eqnarray}
In turn, this relation is the consequence of the well-known Legendre polynomials property,
\begin{eqnarray} 
\int_{-1}^1 P_\ell(x) x^\mu dx=0; \quad \mu<\ell, \label{intPx}
\end{eqnarray}
following from their orthogonality.

\subsection{Discretization of the transform}

Let us introduce the coordinate grid with the step $\Delta r$ in the following way
\begin{eqnarray} 
r_i = (i-1/2) \Delta r; \quad i=1,\ldots,N. \label{r_grid}
\end{eqnarray}
The elements of vector $\bpsi$ of the function values are sampled on the grid as
\begin{eqnarray} 
\psi_i=\psi_\ell(r_i)\sqrt{\Delta r}.
\end{eqnarray}
The elements of the Fourier transform matrix are defined as follows
\begin{eqnarray} 
F_{ni} &=& \frac{1}{\sqrt{(2-\delta_{n0})N}} \frac{e^{ik_nr}-(-1)^\ell e^{-ik_nr_i}}{2i^{\ell+1}(-1)^{\left\lceil\ell/2\right\rceil}} 
\nonumber\\ 
&=&  \frac{1}{\sqrt{(2-\delta_{n0})N}} \times 
\left\{
\begin{array}{ll}
\sin(k_n r_i), &  \text{ even } \ell; \\
\cos(k_n r_i), &  \text{ odd  } \ell.
\end{array}\right.
\label{csFourierBasis}
\end{eqnarray}
Here the momentum grid is 
\begin{eqnarray} 
k_n = n \Delta k; \quad n=p_\ell,\ldots,N_\ell, \label{k_grid}
\end{eqnarray}
where the momentum step is $\Delta k=\pi/\rmax$, the integration interval size is $\rmax=N\Delta r$, and the summation limits are $p_\ell=[1+(-1)^\ell]/2$ and $N_\ell=N+p_\ell-1$, that is $n=0,\ldots,N-1$ for the odd $\ell$ and $n=1,\ldots,N$ for the even ones.

The Fourier transform may be written in the matrix form as 
\begin{eqnarray} 
\mx{f} &=& \mx{F} \bpsi. \label{f_vs_psi}
\end{eqnarray}
It yields as a result the value of the vector of the Fourier expansion coefficients $\mx{f}$ related to the function $\tilde{c}_\ell(k)$ as follows:
\begin{eqnarray}
f_n &=& \tilde{c}_\ell(k_n) \sqrt{w_n}
\end{eqnarray}
where the weights
\begin{eqnarray}
 w_n &=& \left(1-\frac{\delta_{n0}}{2}\right)\Delta k. \label{weight_dk}
\end{eqnarray}
Since the Fourier transform matrix $\mx{F}$ is orthogonal, then under the transform the norm is conserved, that is $\mx{f}^\dag\mx{f}=\bpsi^\dag\bpsi$. The transform \eqref{f_vs_psi} performing through the FFT algorithm requires the number of operations of the order of $N\log_2 N$.

The transform \eqref{c_tildec} conserves the norm according to
\begin{eqnarray}
\int_0^\infty |c_\ell(k)|^2 dk &=& \int_0^\infty |\tilde{c}_\ell(k)|^2 dk \label{conservNormExact}
\end{eqnarray}
The approximation of the integrals in this relation by the trapezoidal rule yields
\begin{eqnarray}
\mx{b}^\dag\mx{b}=\mx{f}^\dag\mx{f}, \label{conservNorm}
\end{eqnarray}
where we introduce the vector $\mx{b}$ composed of the coefficients of the Bessel expansion
\begin{eqnarray}
b_n &=& c_\ell(k_n) \sqrt{w_n}.
\end{eqnarray}

Next, let us write the direct and inverse discrete FtB (DFtB) 
in the following form
\begin{eqnarray}
\mx{b} &=& \mx{T} \mx{f} \label{b_vs_f} \\ 
\mx{f} &=& \mx{T}^{-1} \mx{b}. \label{f_vs_b}
\end{eqnarray}
In order for \eqref{conservNorm} to hold true, the matrix $\mx{T}$ has to be orthogonal, that is
\begin{eqnarray}
\mx{T}^{-1} &=& \mx{T}^T.
\end{eqnarray} 

If we attempt to apply the trapezoidal rule directly to \eqref{c_tildec}, then we would obtain
\begin{eqnarray*}
T_{nm} &=& \delta_{nm} - \frac{P_\ell'(k_m/k_n)}{k_n} \sqrt{w_n}
\end{eqnarray*}
However this technique of the matrix construction does not provide its orthogonality. The reason is that the equation \eqref{intPp_cond} does not hold upon the approximate integration. The employing of the high-order Newton-Cotes rules instead of the trapezoidal rule does not make the situation better. 
In order for \eqref{intPP_cond} to be true, it is necessary for \eqref{intPp_cond} to hold for all the subgrids with the arbitrary nodes number.
The high-order Newton-Cotes rules do not provide high accuracy for an arbitrary subgrid. Therefore the only way to preserve the transform orthogonality appears to be the modification of the integral \eqref{c_tildec} kernel under the proceeding to the numerical integration.

\subsection{Discrete Legendre orthogonal polynomials}

In the context of the summation on grids, the properties analogous to those of the Legendre polynomials are possessed by the so-called discrete Legendre orthogonal polynomials (DLOP) \cite{Neuman1974}.
 DLOP $P_\ell(i,N)$
satisfy the orthogonality property given by
\begin{eqnarray}
 \sum_{i=0}^{N} P_\ell(i,N) P_\mu(i,N) &=& \mathcal{N}(\ell,N) \delta_{\ell\mu} \label{DLOPortho} 
\end{eqnarray}
and also the normalizing condition $P_\ell(0,N)=1$. Here
\begin{eqnarray}
\mathcal{N}(\ell,N) &=&  \frac{(N+\ell+1)^{\underline{\ell+1}}}{(2\ell+1)N^{\underline{\ell}}}, \label{DLOPnormfac}
\end{eqnarray}
where $i^{\underline{j}}$ is $j$-th falling factorial of $i$, $i^{\underline{j}}=i(i-1)\ldots (i-j+1)$.
DLOP might be presented as
\begin{eqnarray}
 P_\ell(i,N) &=& \sum_{j=0}^\ell l(\ell,j) \frac{i^{\underline{j}}}{N^{\underline{j}}},
\end{eqnarray}
where $l(\ell,j)$ are coefficients of the expansion of the shifted Legendre polynomial $P_\ell(1-2x)=\sum_{j=0}^\ell l(\ell,j)x^j$. This means that $P_\ell(i,N)$ can be obtained from $P_\ell(1-2x)$ through the substitution of $i^{\underline{j}}/N^{\underline{j}}$ for $x^j$.
As the grid size increases, DLOP tend to the usual Legendre polynomials according to 
\begin{eqnarray}
 P_\ell(i,N) &=& P_\ell(1-2i/N) +O(N^{-2}).
\end{eqnarray}

Due to the orthogonality condition \eqref{DLOPortho} DLOP possess a property similar to the property \eqref{intPx}, as follows
\begin{eqnarray} 
\sum_{i=0}^{N} P_\ell(i,N) i^s &=& 0; \quad s<\ell. \label{intDPx}
\end{eqnarray}
Making use of this fact one can easily prove (as shown in the Appendix) that for any discrete polynomial $p(i)$ of the order $\mu\leq\ell$ true is the following relation
\begin{eqnarray}
 \sum_{i=0}^{N} P_\ell'(i,N) p(i) w_i(N) &=& (-1)^{\ell} p(N) - p(0) \label{sumDDPp}
\end{eqnarray}
where the weight function coincides with the weights of the trapezoidal integration rule for the grid with a unit step
\begin{eqnarray}
w_i(N) = 1-\frac{\delta_{iN}+\delta_{i0}}{2}.
\end{eqnarray}
Here we define a new discrete polynomial 
\begin{eqnarray}
 P_\ell'(i,N) &=& \frac{2}{1 + P_\ell(-1,N-1)} \nabla [P_\ell](i,N-1)
\end{eqnarray}
which is proportional to the backward difference 
\begin{eqnarray}
 \nabla [P_\ell](i,N-1)=P_\ell(i,N-1)-P_\ell(i-1,N-1) \label{backDiff}
\end{eqnarray}
and hence has the order $\ell-1$. This polynomial tends to the derivative of
the usual Legendre polynomial as
\begin{eqnarray}
 P_\ell'(i,N) &=& \frac{d}{di}P_\ell(1-2i/N) +O(N^{-3}), \label{DPconv}
\end{eqnarray}
that is faster than the DLOP tends to its non-discrete analogue. Therefore we shall further refer to $P_\ell'(i,N)$ as the derivative of the discrete orthogonal Legendre polynomial (DDLOP). It should be mentioned that this term has different meanings throughout the literature \cite{Neuman1974}.

It follows from the relation \eqref{DPconv} that the integral kernel in Eq.\eqref{c_tildec} can be approximated on the grid \eqref{k_grid} by means of DDLOP according to
\begin{eqnarray}
\frac{P_\ell'(k_m/k_n)}{k_n}\Delta k = -P_\ell'(n-m,2n) + O[(\Delta k)^{-3}]. \label{approx_kernel}
\end{eqnarray}

\subsection{Transform matrix}

Let us suppose the transform matrix $\mx{T}$ has the elements as follows:
\begin{eqnarray}
T_{nm} &=& \alpha_n \left[ \delta_{nm} - L_{nm} \right], \label{Tdef}
\end{eqnarray}
where the lower triangular matrix $\mx{L}$ is an approximation of the kernel of the integral \eqref{c_tildec} via \eqref{approx_kernel}, defined by
\begin{eqnarray}
L_{nm} &=& - \theta(n-m) P_\ell'(n-m,2n) \sqrt{1-\frac{\delta_{m0}}{2}}. \label{Ldef}
\end{eqnarray}
Here the Heaviside function $\theta(n-m)$ is specified in \eqref{step_func}, and the additional factor $\sqrt{1-\frac{\delta_{m0}}{2}}$ provides the weight function \eqref{weight_dk} in the product $\mx{T}\mx{f}$. Defined in such a way $T_{nm}$ exist only for $n\geq n_{0\ell}$ where
\begin{eqnarray}
n_{0\ell} = \left\lceil \frac{\ell+1}{2} \right\rceil,
\end{eqnarray}
since there are no DLOP of the order $\ell$ at smaller $n$.

By making use of Eq.\eqref{sumDDPp}, one may show (see the Appendix for details) 
that the rows of the matrix $\mx{I}-\mx{L}$ are mutually orthogonal, that is
\begin{eqnarray}
\sum_{l=p_\ell}^{N_\ell} [\delta_{nl} - L_{nl}] [\delta_{ml} - L_{ml}] &=& \alpha_n^{-2}\delta_{nm};\quad n\geq n_{0\ell}
\label{orthoImL}
\end{eqnarray}
The rows of the matrix $\mx{T}$ are equal to those of the matrix $\mx{I}-\mx{L}$ multiplied by the normalizing constants
\begin{eqnarray}
\alpha_n = \left\{ \sum_l \left[ \delta_{nl} - L_{nl} \right]^2 \right\}^{-1/2}.
\end{eqnarray}
Here the normalizing constants $\alpha_n\to 1$ at $n\to\infty$.

Since the approximation \eqref{approx_kernel} implies the error scaling as $O[(\Delta k)^3]$, the result of the DFtB defined by matrix \eqref{Tdef} differs from the result of the exact FtB by the value $O[(\Delta k)^2]$.

\subsection{Completion of the basis}\label{Sec:Completion}

The Eqs.(\ref{Tdef},\ref{Ldef}) define only $N-n_{0\ell}+1$ rows in the transform matrix. In order to make the basis complete we have to supplement it by extra $n_{0\ell}-1$ vectors that are orthogonal to all other ones.

The DDLOP property Eq.\eqref{sumDDPp} (which Eq.\eqref{orthoImL} follows from) leads to the fact that the basis vectors specified via Eqs.(\ref{Tdef},\ref{Ldef}) are to be orthogonal to any polynomial of the order $s\leq \ell-1$. That is to say, one might construct the extra basis vectors from the 
 polynomials of the order $s\leq \ell-1$. To provide these extra vectors to be orthogonal to not only the basic basis vectors, but to each other as well, we shall choose them as the DLOPs (not the DDLOP!) of the corresponding orders.

So we shall define the extra basis vectors as follows
\begin{eqnarray}
T_{nm} &=& \alpha_n P_{l(n)}(N_\ell-m,2N_\ell) \sqrt{1-\frac{\delta_{m0}}{2}}; \, n\in\left[p_\ell, n_{0\ell}-1 \right], \label{Tadd}
\end{eqnarray}
where the polynomial order is
\begin{eqnarray}
 l(n) &=& 2n-p_\ell,
\end{eqnarray}
and the normalizing constant is $\alpha_n=\sqrt{2/\mathcal{N}[l(n),2N_\ell]}$, where $\mathcal{N}$ is given in Eq.\eqref{DLOPnormfac}.

Now let us elucidate the physical meaning of the extra basis vectors \eqref{Tadd}. For this purpose we shall derive their appearance in coordinate representation. First, let us begin with the expression for the coefficient
\begin{eqnarray}
 b_n &=& \sum_{m=p_\ell}^{N_\ell} \alpha_n P_{l(n)}(N_\ell-m,2N_\ell) \sqrt{1-\frac{\delta_{m0}}{2}} f_m \label{lown_b_vs_f} 
\end{eqnarray}
and proceed making use of the relation between the DLOP and the usual Legendre polynomial
\begin{eqnarray*}
 P_\ell(N-m,2N) &\approx& P_\ell\left(m/N\right).
\end{eqnarray*}
By changing the summation to integration we get
\begin{eqnarray*}
 b_n &\approx& \alpha_n \int_0^{k_{\max}} P_{l(n)}(q/k_N) \tilde{c}_\ell(q) dq,
\end{eqnarray*}
where $k_N=N\Delta k$. 
Finally, substitute here \eqref{tilde_c_def} and use Eq.\eqref{BesselExp} to obtain
\begin{eqnarray}
b_n &\approx& \frac{k_N\alpha_n}{i^{\ell+1-l(n)}}\sqrt{\frac{2}{\pi}} \int_0^\infty j_{l(n)}(k_Nr) \psi_\ell(r) dr \label{addFuncProj}
\end{eqnarray}
Thus, $b_n$ at $n<n_{0\ell}$ are the coefficients of the expansion in the functions $j_{l(n)}(k_Nr)$ (while $b_n$ at $n\geq n_{0\ell}$ are the coefficients of the $\psi_\ell(r)$ expansion in $\chi_{\ell}(k_nr)$). At $r\to 0$ asymptotics are $\chi_{\ell}(kr) \sim r^{\ell+1}$ and $j_{l(n)}(k_Nr)\to r^{l(n)}$.
Since $l(n)\leq\ell-1$, the additional basis vectors represent the high-momentum wave-function components converging to zero slower than $r^{\ell+1}$. That is to say, the additional vectors emerge to be a linear combination of the non-regular SBFs (SBFs of the second kind).

If $\psi_\ell(r)$ appears to be a result of the spherical harmonics expansion of a 3D function having the continuous derivatives up to order $\ell+1$, then $\psi_\ell(r\to 0)\sim r^{\ell+1}$, hence it should be $b_n=0$ at $n< n_{0\ell}$. In the case when the function possesses singularities such components become non-zero and should be considered as having the energy larger than that of the component specified by the coefficient $b_N$.

Thus, our DSBT does not yield coefficients for the small momenta $k_n, n< n_{0\ell}$. The reason for this is that at such $k_n$ there exist no SBF of the first kind satisfying the boundary conditions (see also the Section \ref{Sec:Testing}).

\subsection{Fast multiplication by transform matrix}

The fast transform might be accomplished by means of the technique proposed in the work \cite{Toyoda2010}, except that we shall use DDLOP instead of the Legendre polynomials. 

First, let us expand DDLOP over the powers of the variable $m$ in the following way:
\begin{eqnarray}
P_\ell'(n-m,2n) = -\sum_{\nu=0}^{\ell-1} \xi_{\ell\nu}(n) m^{\nu}. \label{DP_expan}
\end{eqnarray}
Upon substituting this expansion into \eqref{b_vs_f} one obtains
\begin{eqnarray}
b_n &=& \alpha_n f_n - \alpha_n \sum_{\nu=0}^{\ell-1} \xi_{\ell\nu}(n) s_{\nu n}. \label{fast_b_f}
\end{eqnarray}
Here we introduce the notation
\begin{eqnarray}
s_{\nu n} &=& \frac{\delta_{\nu 0}}{\sqrt{2}}f_0 + \sum_{m=1}^{n-1} m^{\nu} f_m + \frac{1}{2} n^{\nu} f_n. \label{sum_in_b_f} 
\end{eqnarray}
The sum \eqref{sum_in_b_f} can be evaluated through the recurrence relation 
\begin{eqnarray}
s_{\nu n} &=& s_{\nu, n-1} + \frac{1}{2}\left[n^{\nu} f_n + (n-1)^{\nu} f_{n-1}\right]. \label{rec_sum_b_f}
\end{eqnarray}
Thereby the sums calculation for all the $n\in [n_{0\ell},N]$ requires as less as $N$ operations. In total one needs to evaluate $\left\lceil\ell/2\right\rceil\times (N+1-n_{0\ell})$ sums ($\left\lceil\ell/2\right\rceil$ appearing because of the fact that the half of the coefficients in \eqref{DP_expan} are zero due to the DDLOP parity), then to perform the summation over $\nu$ in \eqref{fast_b_f}, resulting in the altogether operations number scaling as $\ell N$.

At $n<n_{0\ell}$ the coefficients $b_n$ are to be computed according to \eqref{lown_b_vs_f}. That means that for every $n<n_{0\ell}$ one has to calculate the vectors scalar product that requires extra $O(\ell N)$ operations, that is the operations number scaling is the same as for the $b_n$ evaluating for $n\geq n_{0\ell}$ via (\ref{rec_sum_b_f},\ref{fast_b_f}). In sum, to accomplish the transform \eqref{b_vs_f} one needs to perform $O(\ell N)$ operations.

Now let us consider the inverse transform. The substitution of \eqref{DP_expan} into \eqref{f_vs_b} yields for $m\geq n_{0\ell}$ the following:
\begin{eqnarray}
f_m &=& \sum_{n=p_\ell}^{n_{0\ell}-1} T_{nm} b_n + \alpha_m b_m - \sum_{\nu=0}^{\ell-1} m^{\nu} \tilde{s}_{\nu m}, \label{fast_f_b}
\end{eqnarray}
where
\begin{eqnarray}
\tilde{s}_{\nu m} &=& \sum_{n=m+1}^{N_\ell} \alpha_n \xi_{\ell\nu}(n) b_n + \frac{1}{2}\alpha_m \xi_{\ell\nu}(m) b_m. \label{sum_in_f_b} 
\end{eqnarray}
Next, for $m<n_{0\ell}$ we obtain
\begin{eqnarray}
f_m &=& \sum_{n=p_\ell}^{n_{0\ell}-1} T_{nm} b_n - \sum_{\nu=0}^{\ell-1} \sqrt{1-\frac{\delta_{m0}}{2}} m^{\nu} \tilde{s}_{\nu, n_{0\ell}+1/2} \label{fast_f_b_lowm} 
\end{eqnarray}
where
\begin{eqnarray}
\tilde{s}_{\nu, n_{0\ell}+1/2} &=& \tilde{s}_{\nu n_{0\ell}}+\frac{1}{2}\alpha_{n_{0\ell}} \xi_{\ell\nu}(n_{0\ell}) b_{n_{0\ell}}
\end{eqnarray}
The sum \eqref{sum_in_f_b} may be evaluated according to the recurrence relation 
\begin{eqnarray}
\tilde{s}_{\nu m} &=& \tilde{s}_{\nu, m+1} + \frac{1}{2}\left[\alpha_m \xi_{\ell\nu}(m) b_m + \alpha_{m+1} \xi_{\ell\nu}(m+1) b_{m+1}\right]. \label{rec_sum_f_b}
\end{eqnarray}
Thus, the inverse transform \eqref{f_vs_b} performing requires the number of operations $O(\ell N)$, just as the direct transform does.

The fast Fourier transform \eqref{f_vs_psi} operations number scales as $N\log_2 N$. That is the DSBT in total requires $O(N\log_2 N)+O(\ell N)$ operations. At large $N$ the FFT strongly over-demands the DFtB, hence the overall DSBT algorithm processing time appears to be defined by the FFT processing time.

\section{Numerical test of the method} \label{Sec:Testing}

The method convergence was examined on three grids \eqref{r_grid} with the same space step $\Delta r=0.4$ and the various values of the space region $\rmax=N\Delta r=$ 51.2, 102.4, and 204.8. That is, the grids possessed the same maximal momentum $k_N=\Delta k N=\pi/\Delta r$ and various momentum steps $\Delta k=\pi/\rmax$.

\begin{figure}[ht]
\includegraphics[angle=-90,width=0.45\columnwidth]{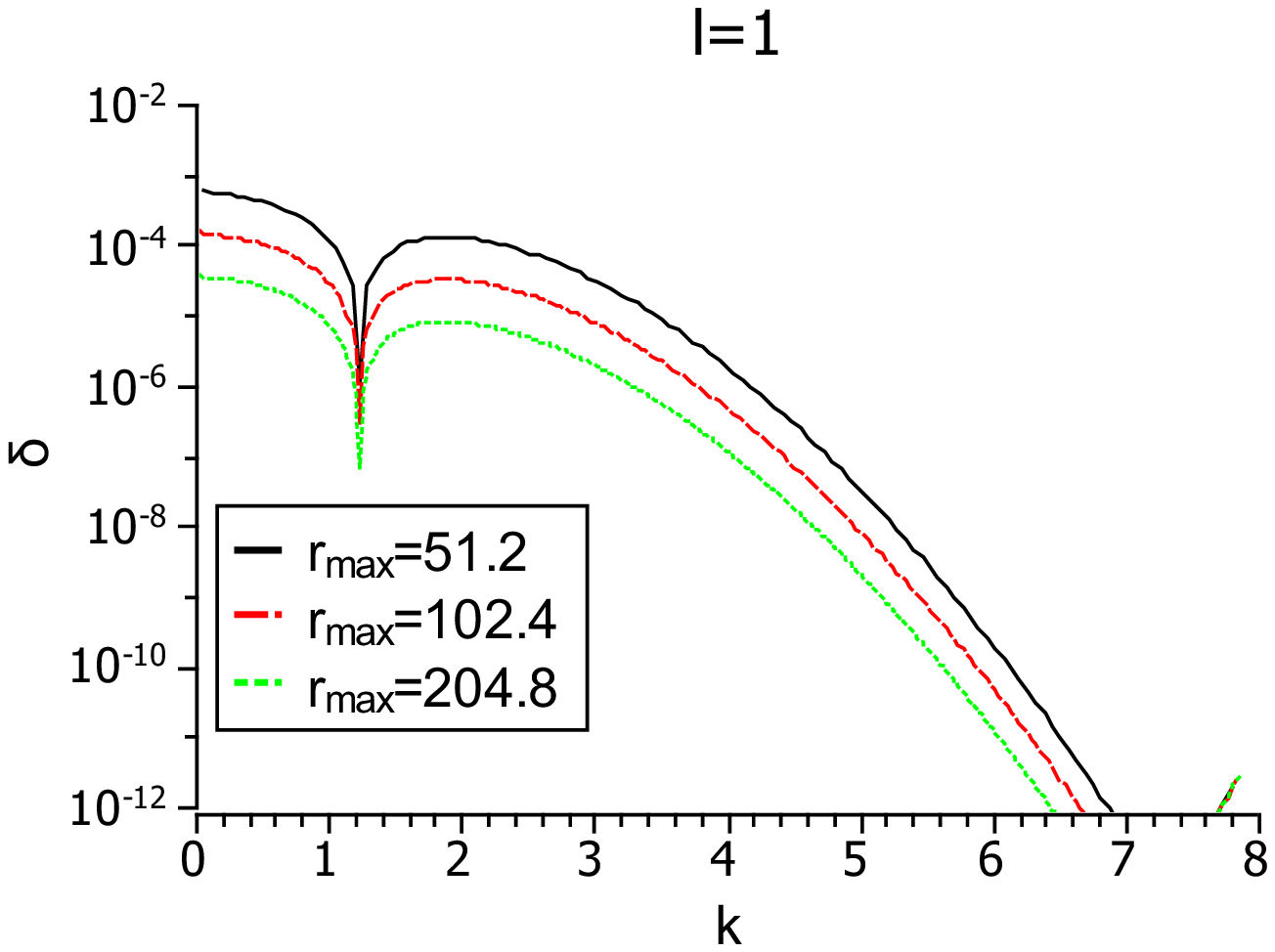}
\includegraphics[angle=-90,width=0.45\columnwidth]{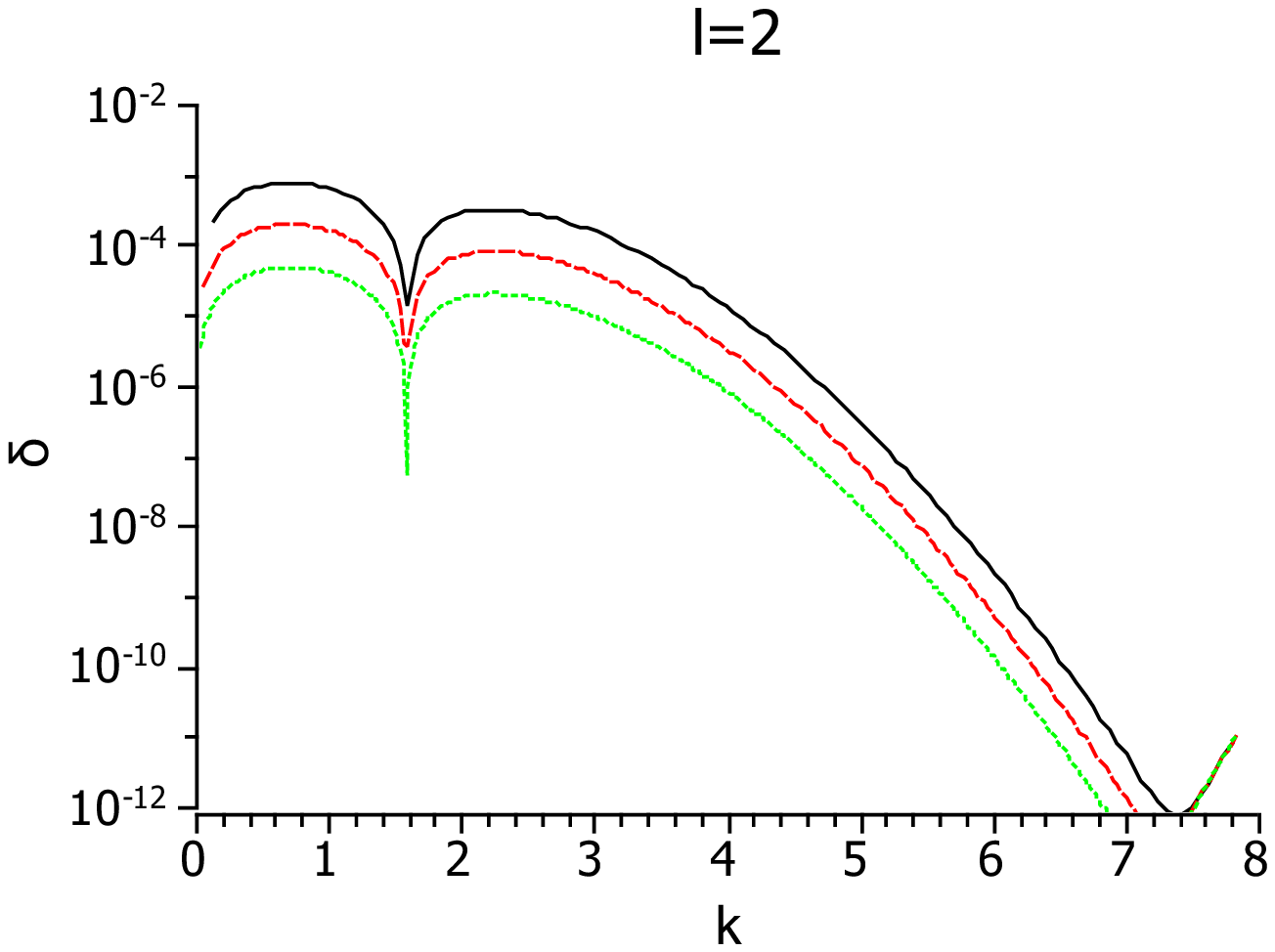}
\\
\includegraphics[angle=-90,width=0.45\columnwidth]{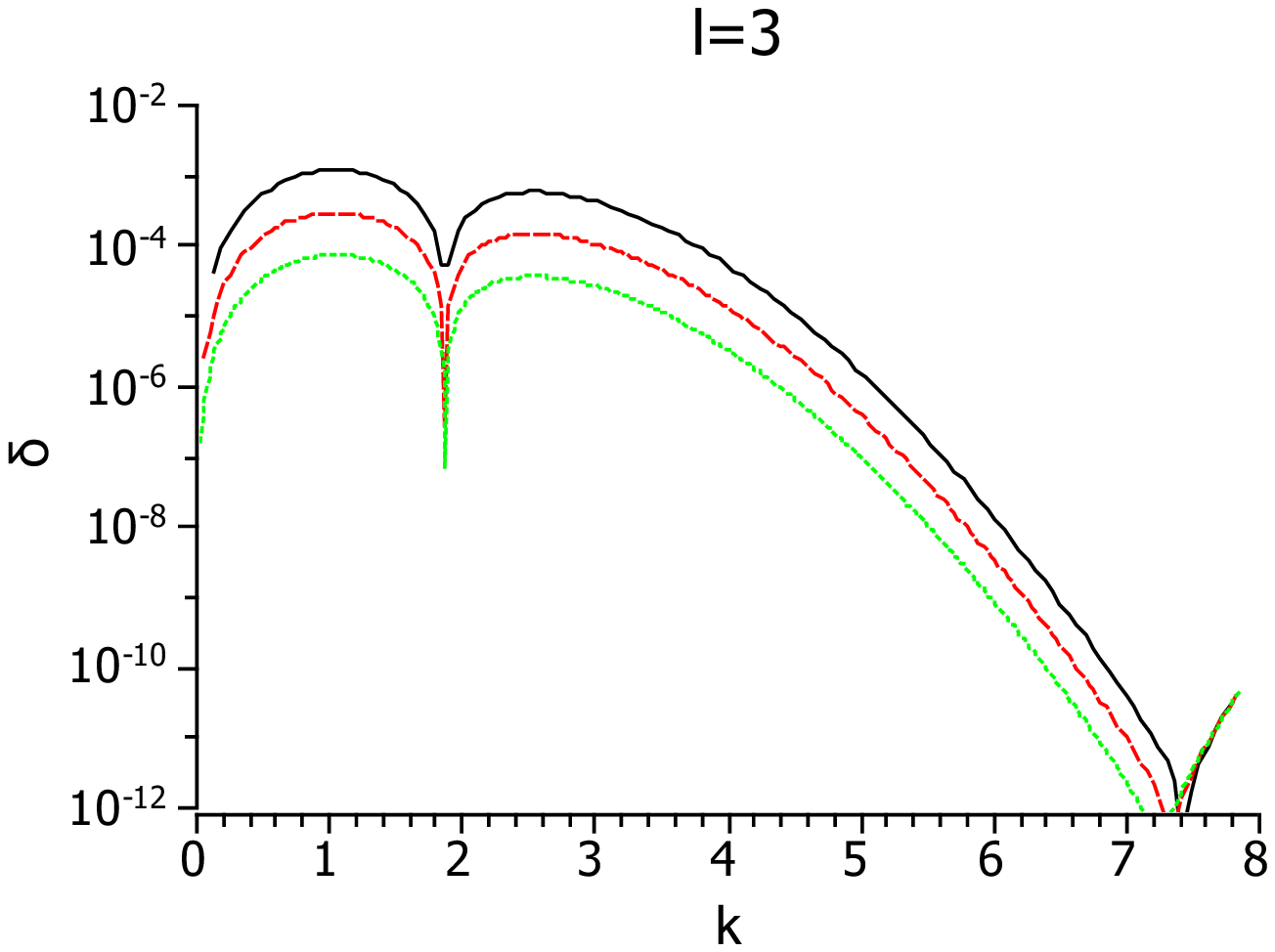}
\includegraphics[angle=-90,width=0.45\columnwidth]{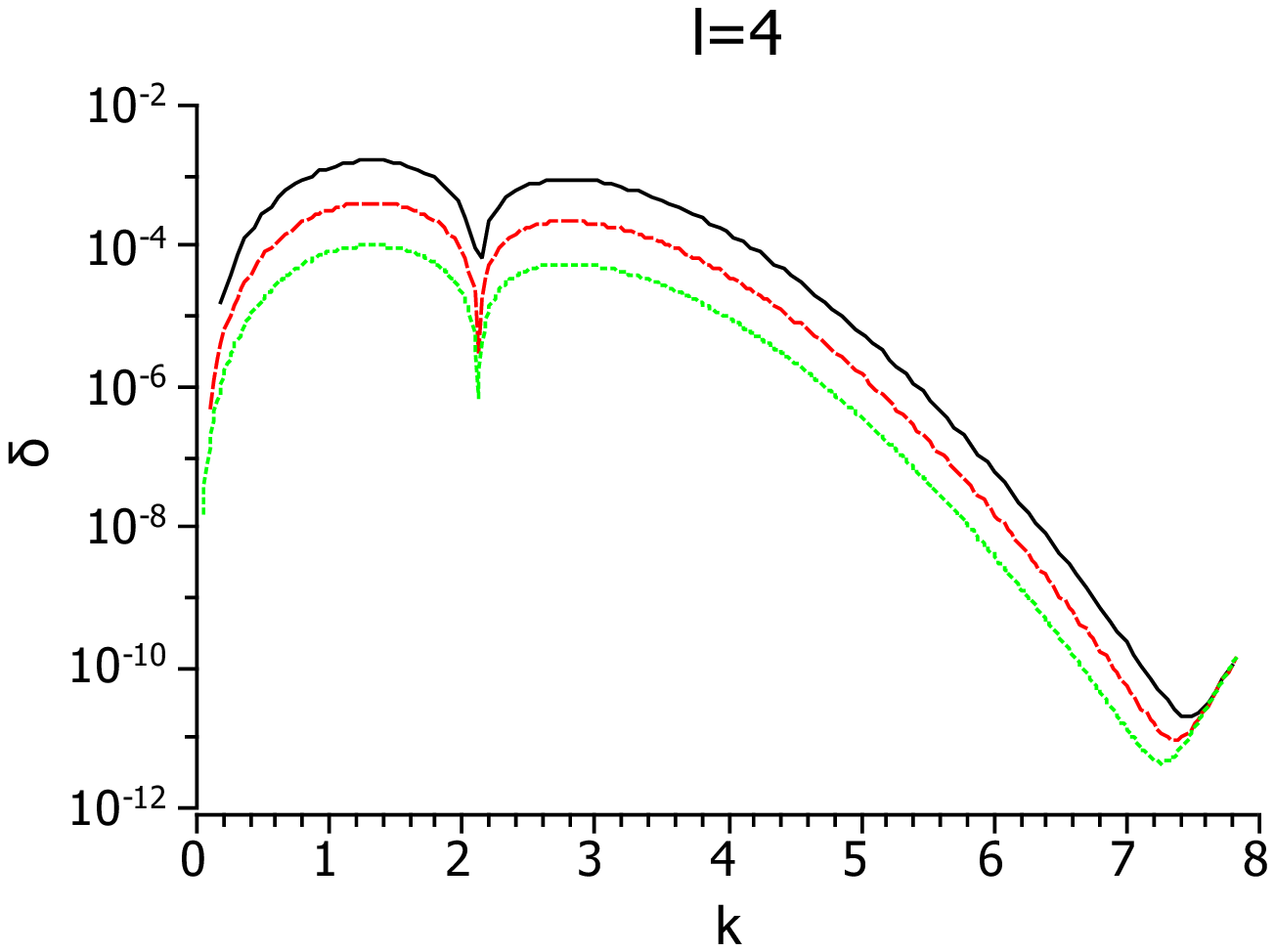}
\caption{Deviation of the Gaussian atomic orbitals expansion from exact for $\rmax$=51.2 (solid lines), $\rmax$=102.4 (dashed lines) and $\rmax$=204.8 (dotted lines)}
\label{fig:Proj}
\end{figure}
Let us begin with the check of the convergence of our transform for smooth functions, which are commonly used in atomic physics namely Gaussian atomic orbital functions
\begin{eqnarray}
\psi_\ell(r) = A_\ell r^{\ell+1}\exp(-r^2/2).
\end{eqnarray}
The Fig.\ref{fig:Proj} shows the absolute value $\delta=|c_{n\ell}-c_{\ell}(k_n)|$ of the difference between the exact SBT result $c_{\ell}(k)$ (obtained by the numerical evaluation of the integral in Eq.\eqref{BesselProj}) and the DSBT result $c_{n\ell} = \Delta k^{-1/2} [\mx{T}\mx{F}\bpsi]_{n}$.
It is seen that upon the $\rmax$ increasing the error decreases as $1/\rmax^{2}$ (or, equivalently, as $\Delta k^2$). It coincides with the convergence rate expected from theory. 

\begin{figure}[ht]
\includegraphics[angle=-90,width=0.45\columnwidth]{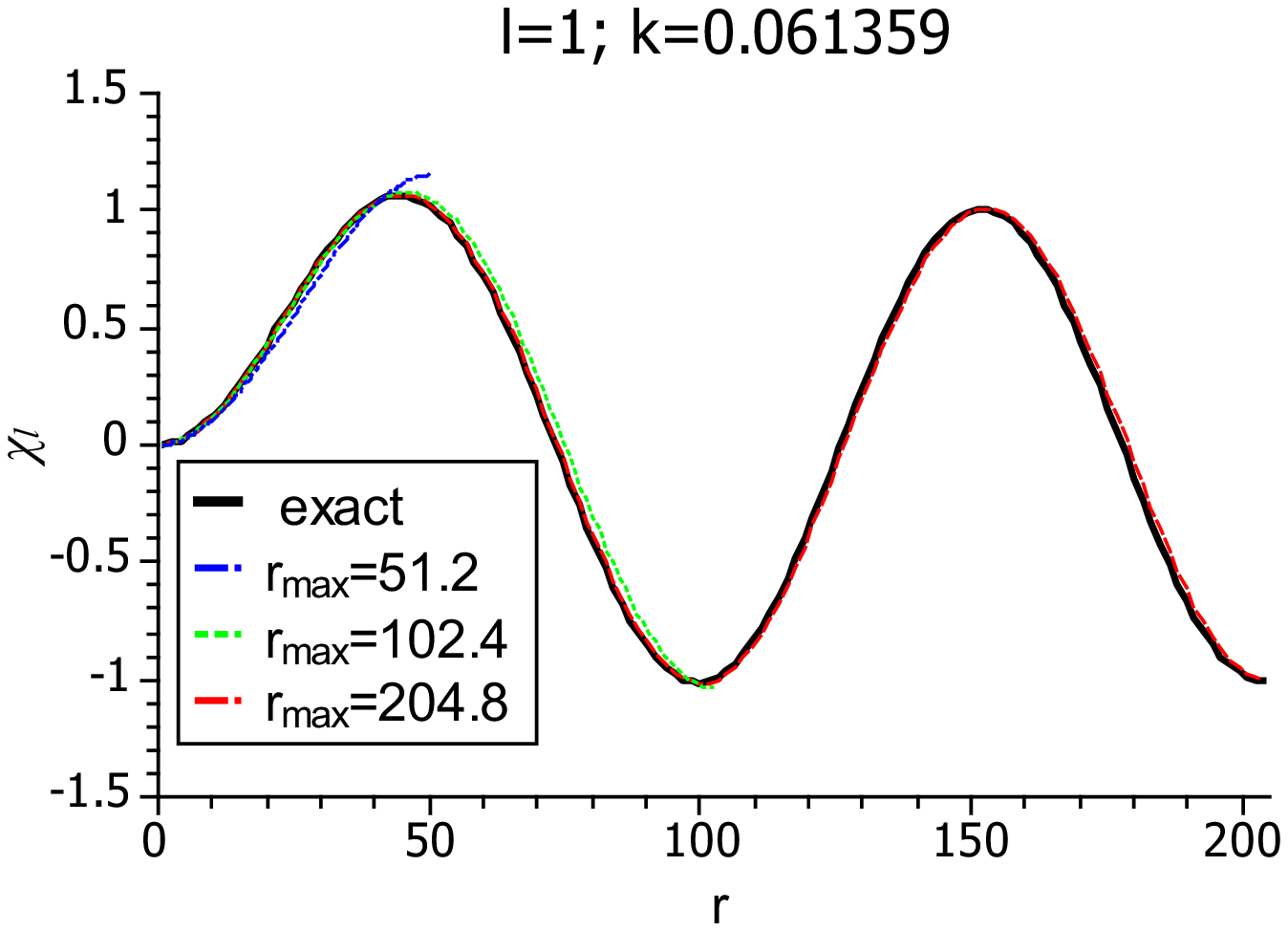}
\includegraphics[angle=-90,width=0.45\columnwidth]{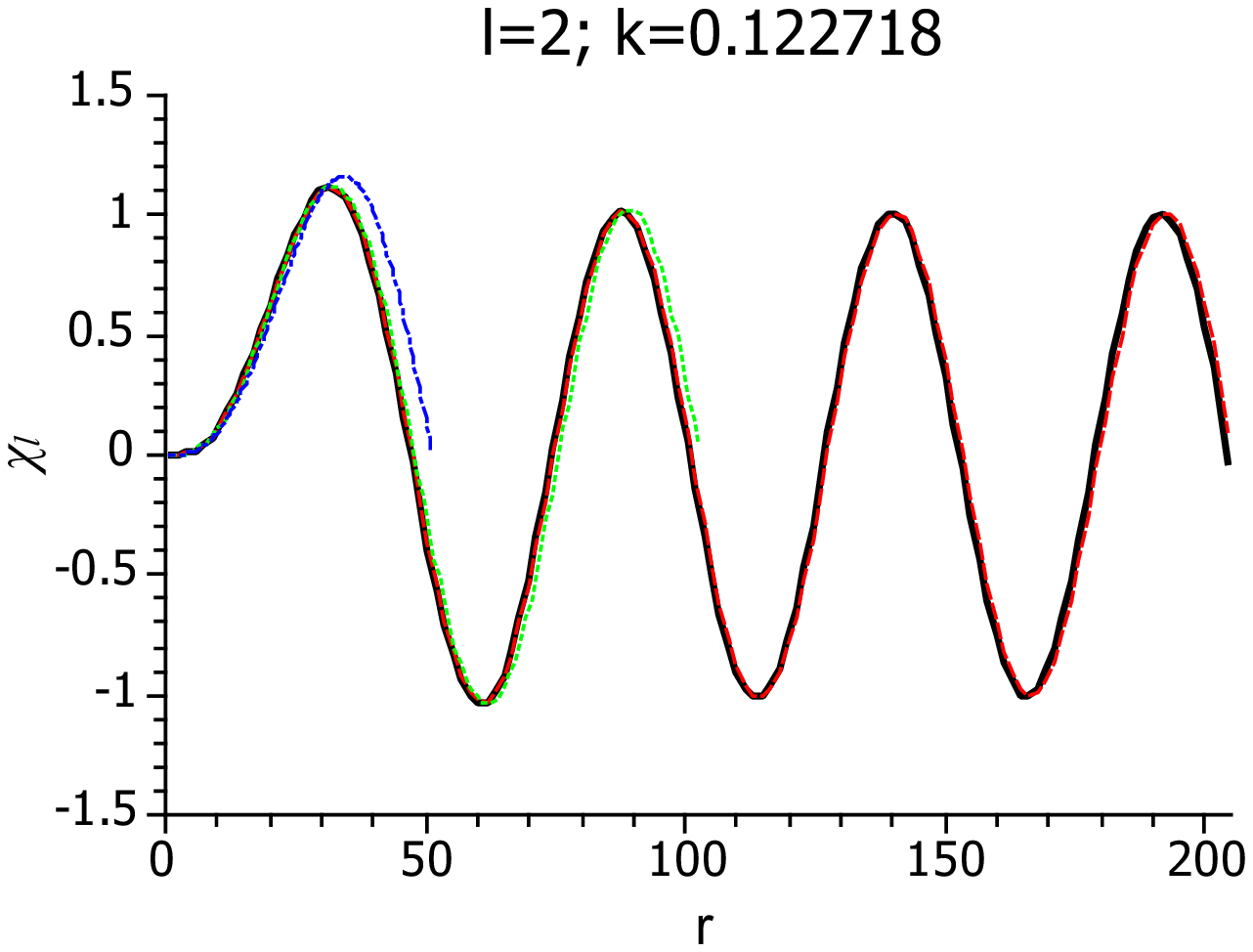}
\\
\includegraphics[angle=-90,width=0.45\columnwidth]{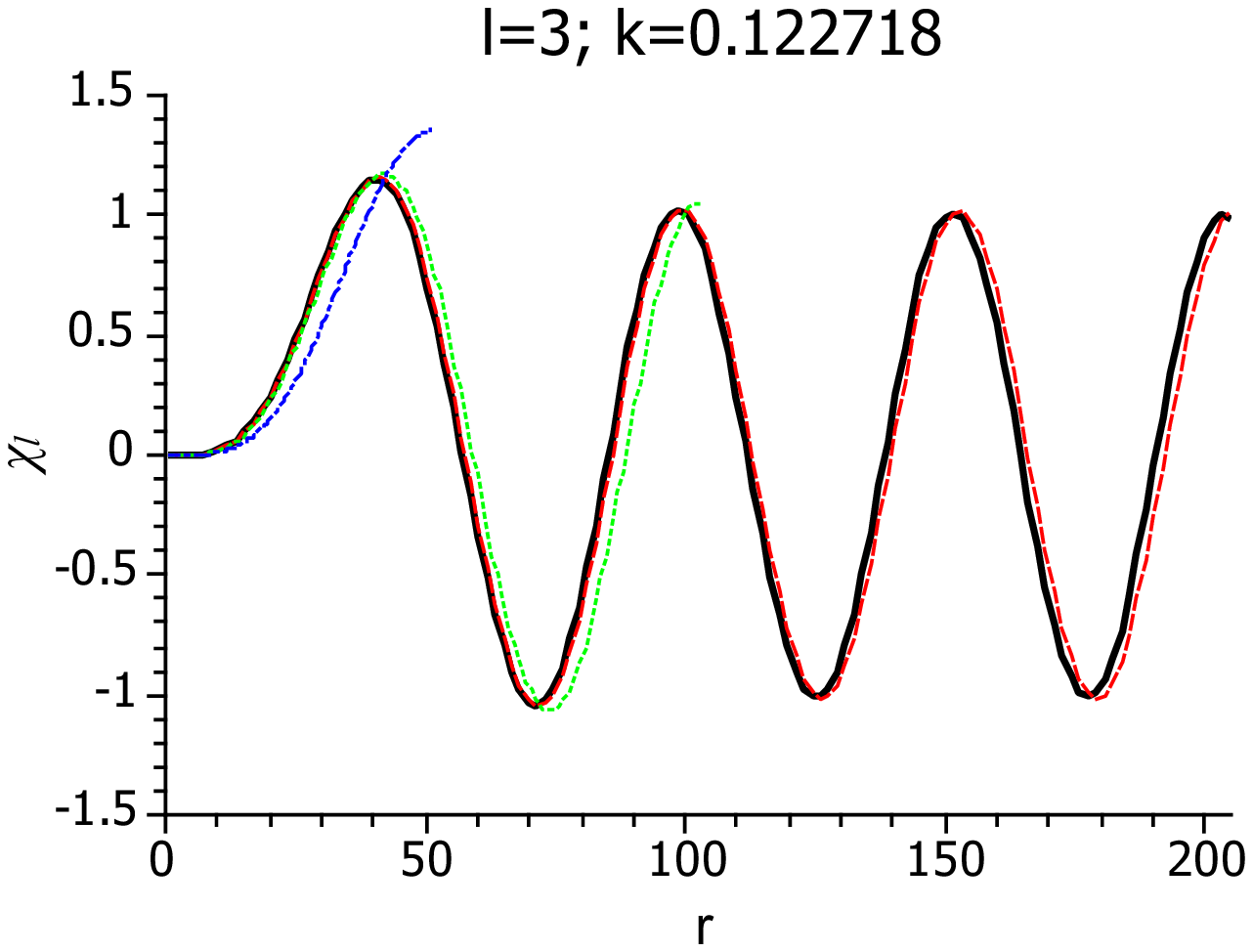}
\includegraphics[angle=-90,width=0.45\columnwidth]{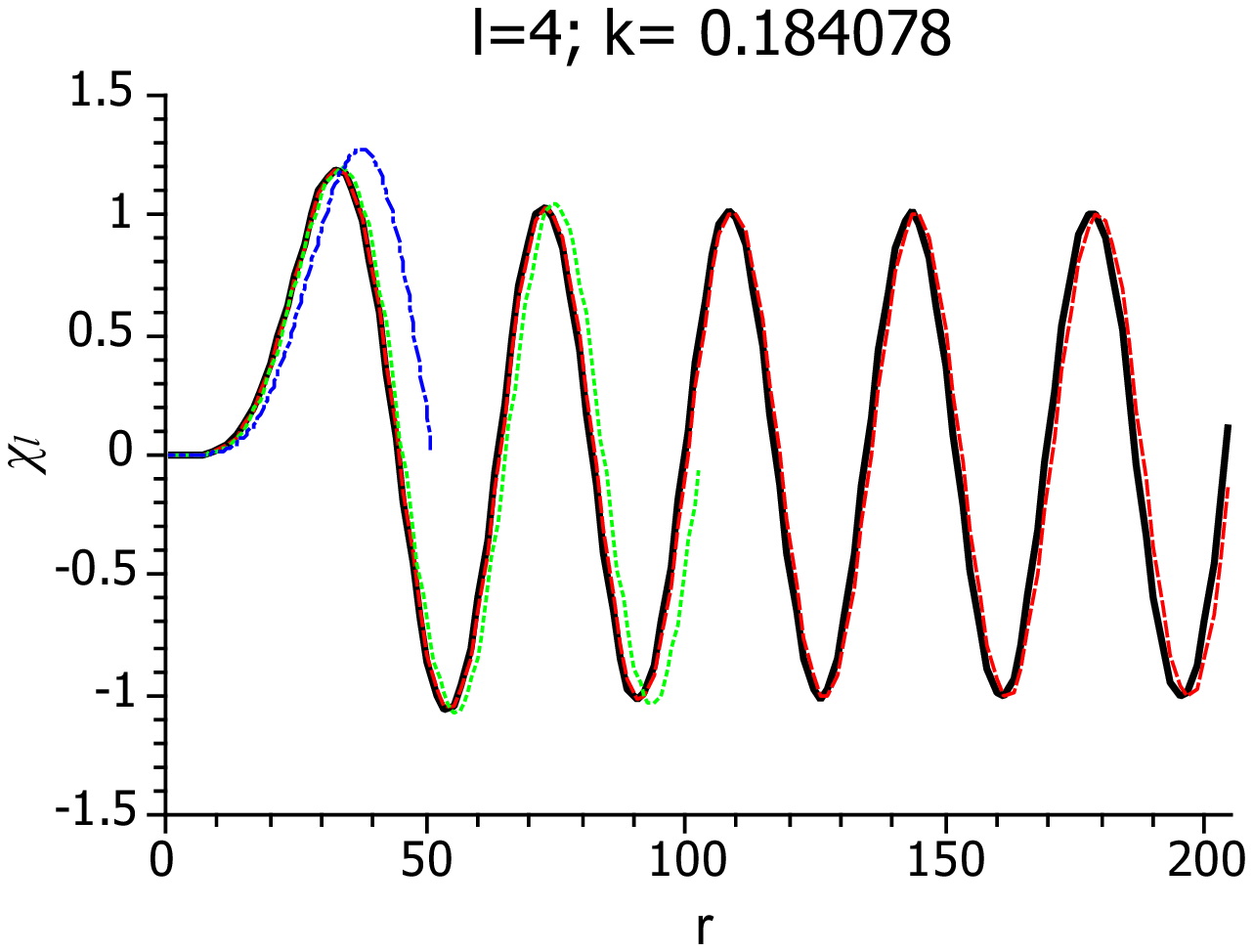}
\caption{Convergence of the DSBT basis functions to SBFs with decreasing of $\Delta k$.}
\label{fig:bfuncs}
\end{figure}
In order to reveal the main error source, let us perform the comparison of the exact SBFs with 
the DSBT basis functions in the coordinate representation, which correspond to the transform basis vectors.
For this purpose we apply the inverse transform to the vector composed of coefficients $b_m=(\rmax/2)^{1/2} w_n^{-1/2} \delta_{mn(k)}$, where $n(k)=k/\Delta k$ and $k$ is the fixed
momentum. The result function may be written as follows:
\begin{eqnarray} 
\tilde{\chi}_{n\ell}(r_i) = (\rmax/2)^{1/2} w_n^{-1/2} [\mx{F}^\dag \mx{T}^T]_{i,n}. 
\end{eqnarray}
The factor $(\rmax/2)^{1/2}$ is introduced in order to provide the convergence $\tilde{\chi}_{n\ell}(r_i) \to \chi_{\ell}(k_nr_i)$ at $\rmax\to \infty$. The results are presented at the Fig.\ref{fig:bfuncs}. For every $\ell$ we chose the momentum $k$ in such a way that for the maximal $\Delta k$ grid a basis vector number $n(k)=n_{0\ell}$. The $n_{0\ell}$-th function has the nodes number minimal among the functions satisfying the boundary conditions at $r=\rmax$ and the asymptotics $\sim r^{\ell+1}$ at $r=0$. Upon the step $\Delta k$ decreasing at the fixed $k$ the number $n(k)$ grows and the DSBT basis function in its domain of definition becomes closer to the SBF.

\begin{figure}[ht]
\begin{center}
\includegraphics[angle=-90,width=0.45\columnwidth]{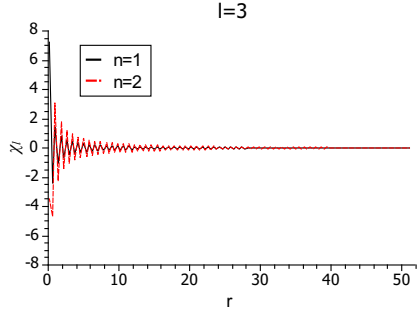}
\caption{Extra basis functions for $\ell=3$.}
\end{center}
\label{fig:afuncs}
\end{figure}
The Fig.\ref{fig:afuncs} demonstrates the DSBT basis functions at $n<n_{0\ell}$. One can see these functions behavior is just as expected from \eqref{addFuncProj}, that is they are high-frequency, decrease with distance and have asymptotic behavior different from $\sim r^{\ell+1}$ at $r=0$.

\begin{figure}[ht]
\includegraphics[angle=-90,width=0.45\columnwidth]{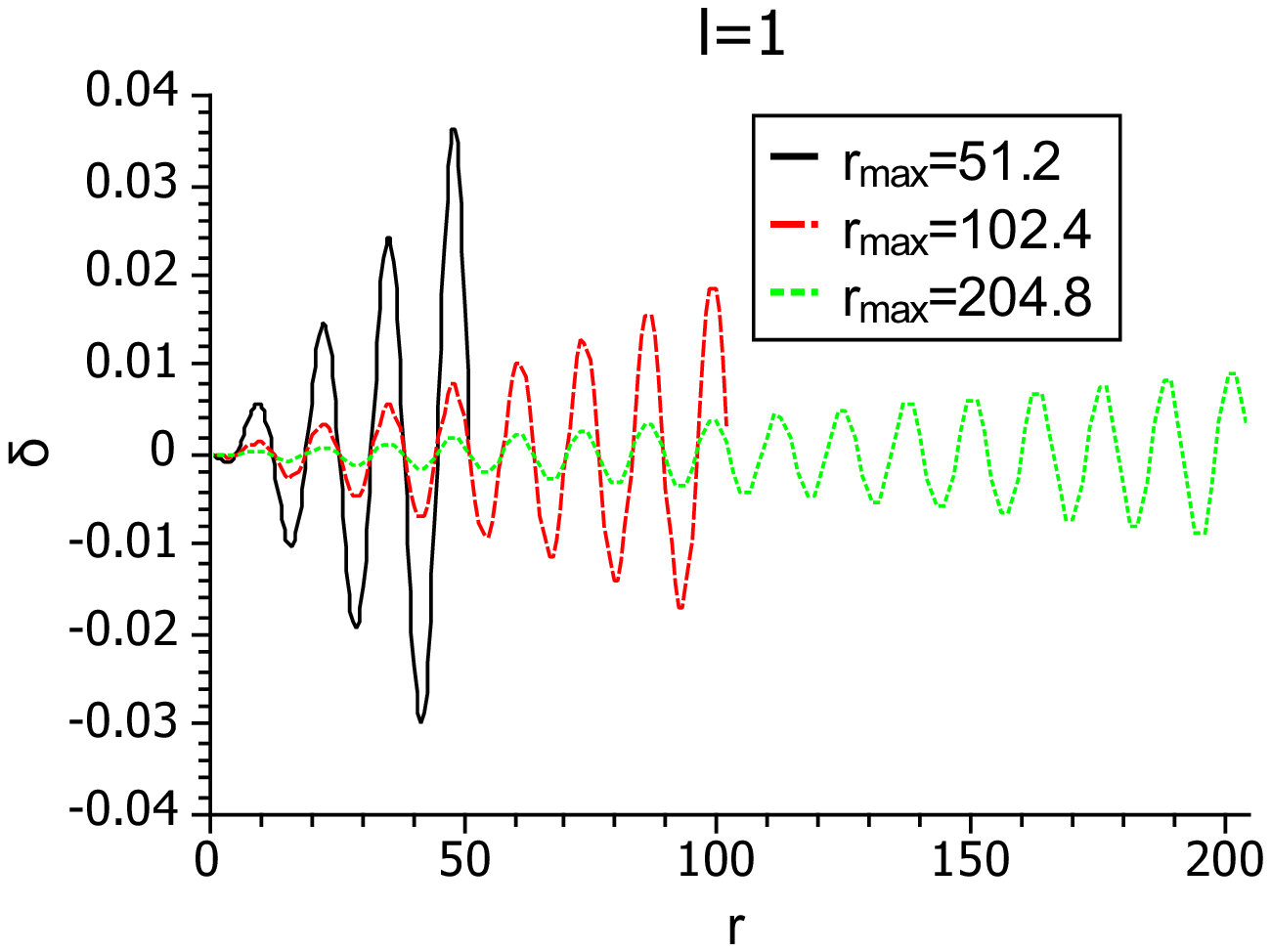}
\includegraphics[angle=-90,width=0.45\columnwidth]{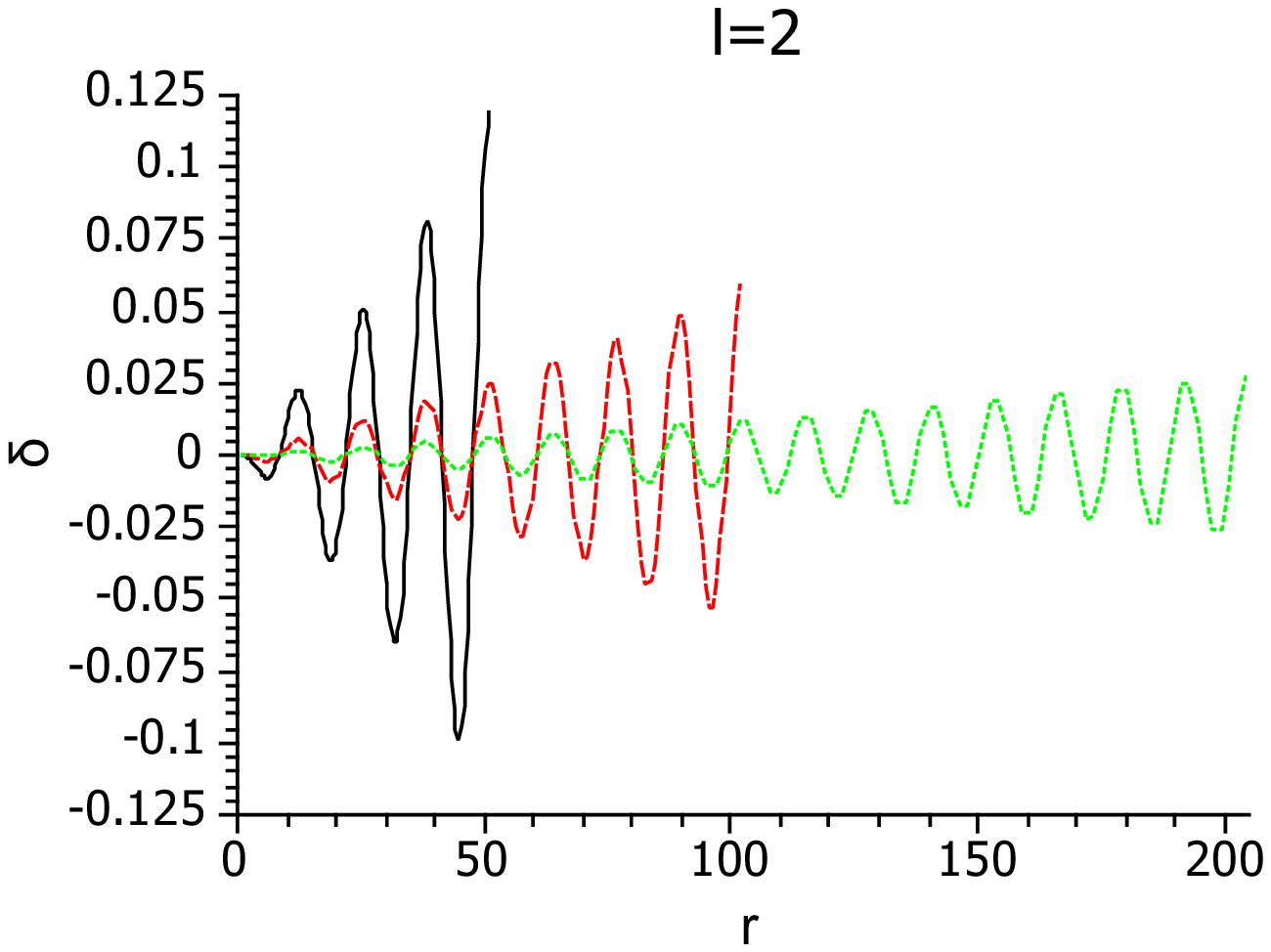}
\\
\includegraphics[angle=-90,width=0.45\columnwidth]{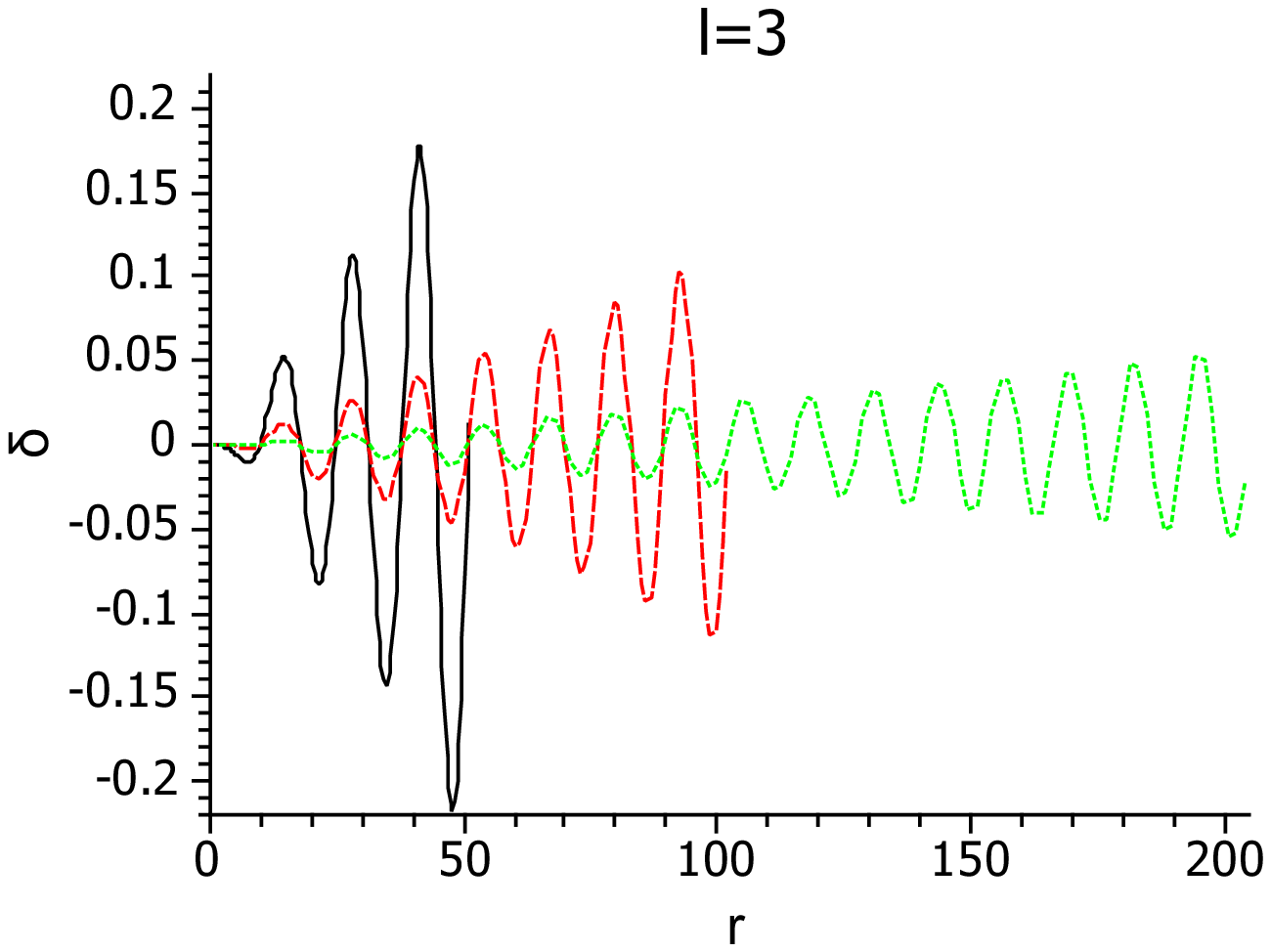}
\includegraphics[angle=-90,width=0.45\columnwidth]{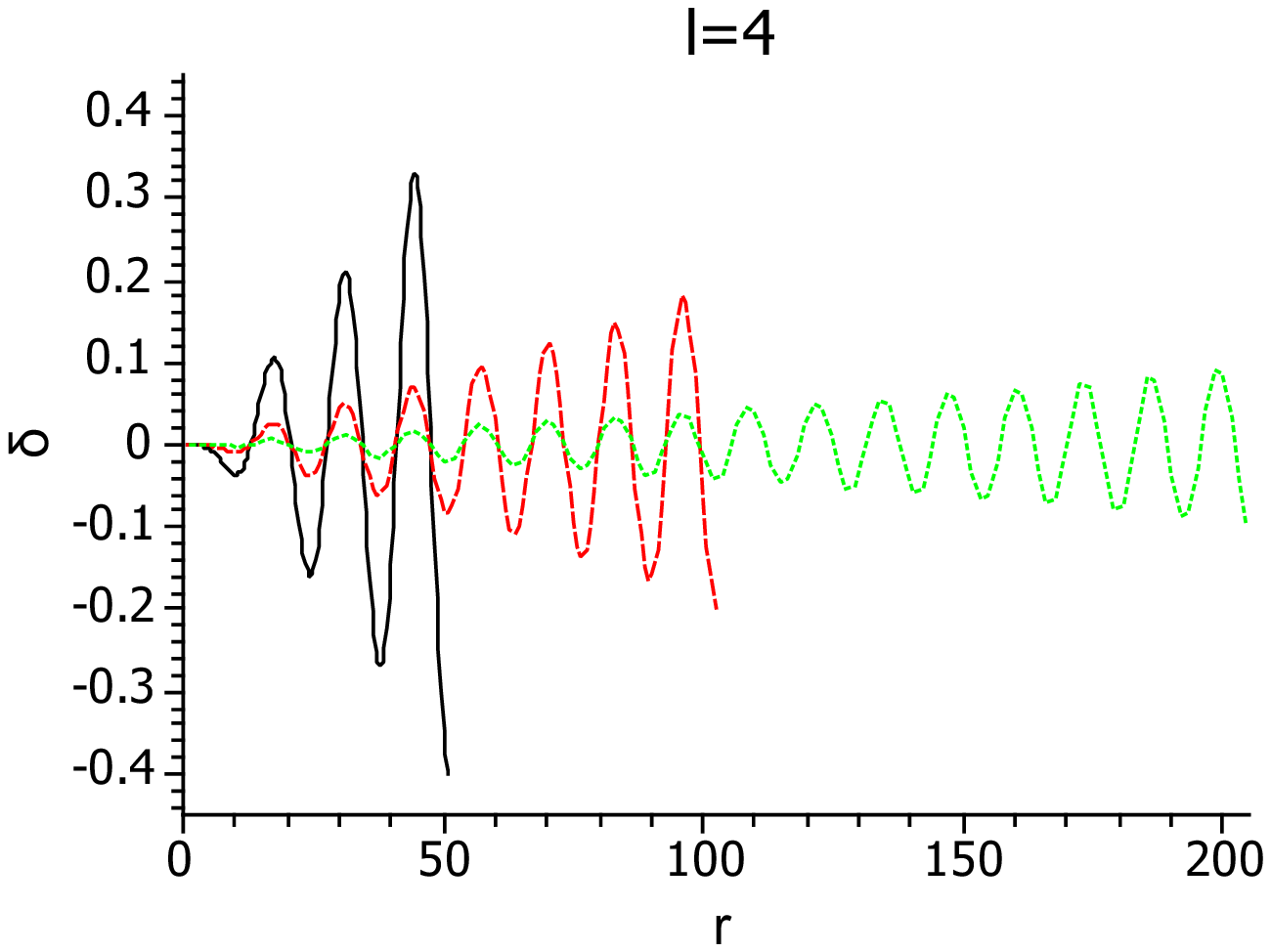}
\caption{Difference between the basis functions and SBFs $\delta=|\tilde{\chi}_{n(k)\ell}(r)-\chi_{\ell}(kr)|$ for the fixed $k$ = 0.490873843.}
\label{fig:berr}
\end{figure}
In order to demonstrate the main error source, we plot the dependence of $\delta=|\tilde{\chi}_{n(k)\ell}(r)-\chi_{\ell}(kr)|$ on radius (Fig.\ref{fig:berr}). At that we take the momentum value such that $n(k)=9$ even at the smallest grid. It is apparent that the difference grows with distance. The reason for the difference to be maximal at $r=r_N$ may be understood if one recalls that SBFs asymptotically equivalent to 
\begin{eqnarray} 
\chi_{\ell}(kr) \simeq \sin[kr-\ell\pi/2-l(l+1)/2kr]; \quad r\to\infty. 
\end{eqnarray}
The DSBT basis function satisfies the same boundary conditions as the functions constituting it (sines or cosines) do
\begin{eqnarray} 
\tilde{\chi}_{n(k)\ell}(r) \simeq \sin(kr-\ell\pi/2); \quad r\to r_N. 
\end{eqnarray}
That means that there does exist a phase shift $l(l+1)/2kr_N$ near the outer boundary between the approximate and exact functions.

The fact of the amplitude of the difference between the SBFs and DSBT basis functions growing roughly linearly with $r$ increasing (which is seen from Fig.\ref{fig:berr}) indicates that the phase shift between these functions grows linearly with $r$ as well. Such phase shift behavior may be interpreted as a consequence of the approximate and exact functions wavenumbers being distinct. Therefore the DSBT basis function has to be closer to the SBF at such momentum $k_{n\ell}$ that provides the exact SBF coincidence with the DSBT basis function on the grid boundary. This condition might be written mathematically as
\begin{eqnarray}
\begin{array}{ll}
\chi_\ell(k_{n\ell}r_{max})=0, &  \text{ even } \ell; \\
\chi_\ell'(k_{n\ell}r_{max})=0, &  \text{ odd  } \ell.
\end{array} \label{kcorr}
\end{eqnarray}
From the phase shift between the DSBT basis functions and exact SBFs near the outer boundary one can obtain the approximate expression
\begin{eqnarray}
k_{n\ell}\simeq k_n - \frac{\ell(\ell+1)}{2\pi^2k_n}\Delta k^2.
\end{eqnarray}

\begin{figure}[ht]
\includegraphics[angle=-90,width=0.45\columnwidth]{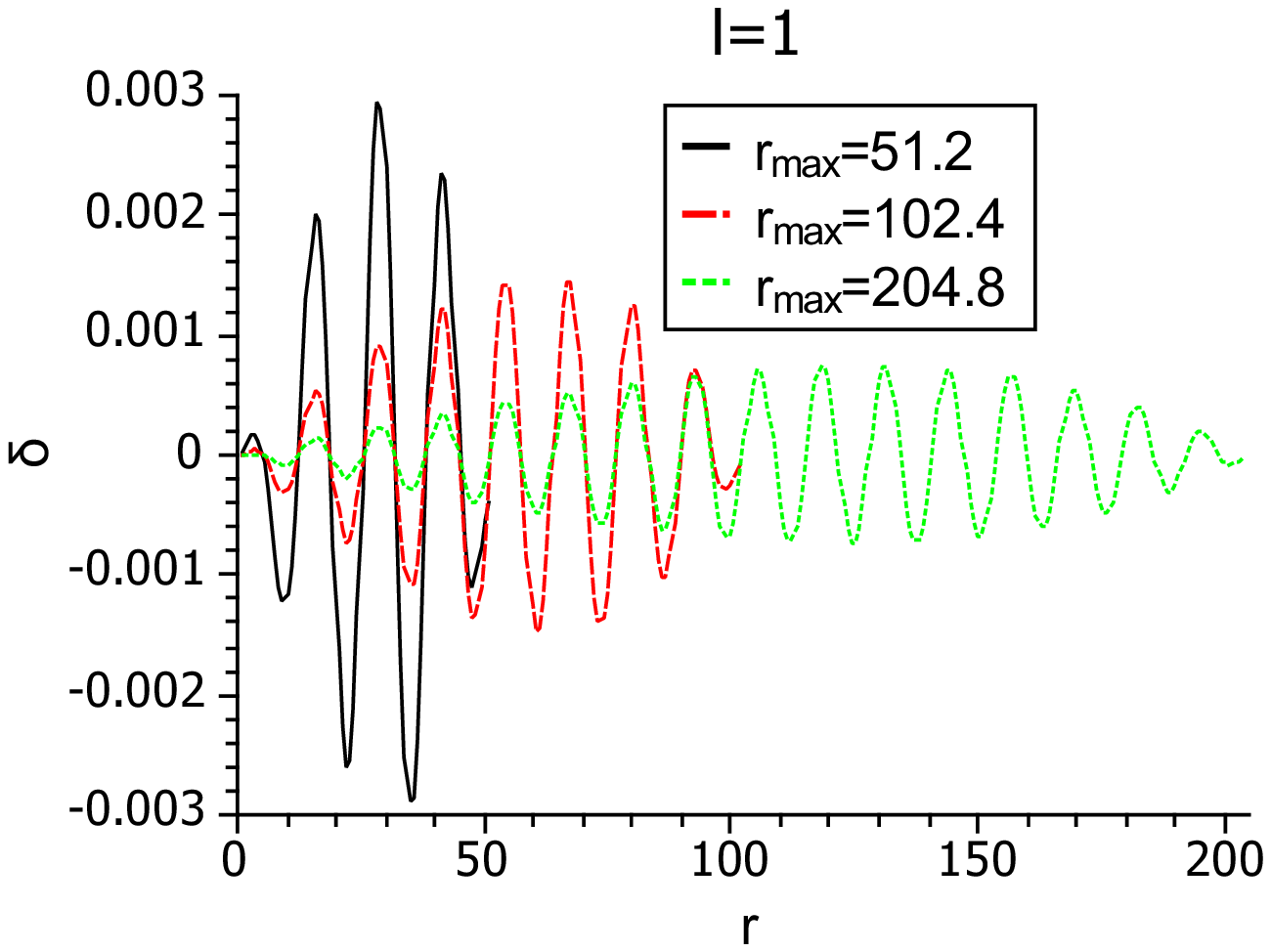}
\includegraphics[angle=-90,width=0.45\columnwidth]{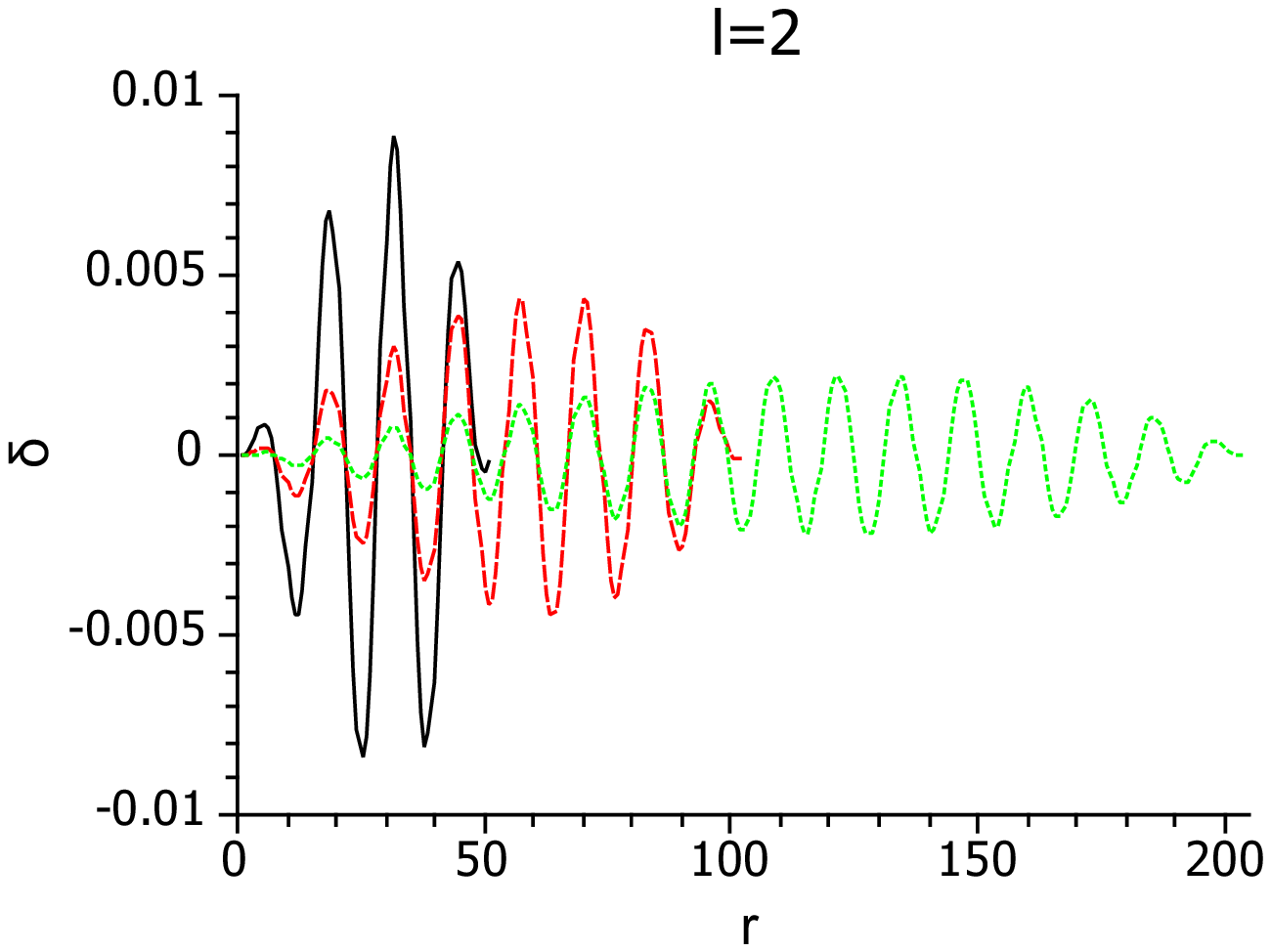}
\\
\includegraphics[angle=-90,width=0.45\columnwidth]{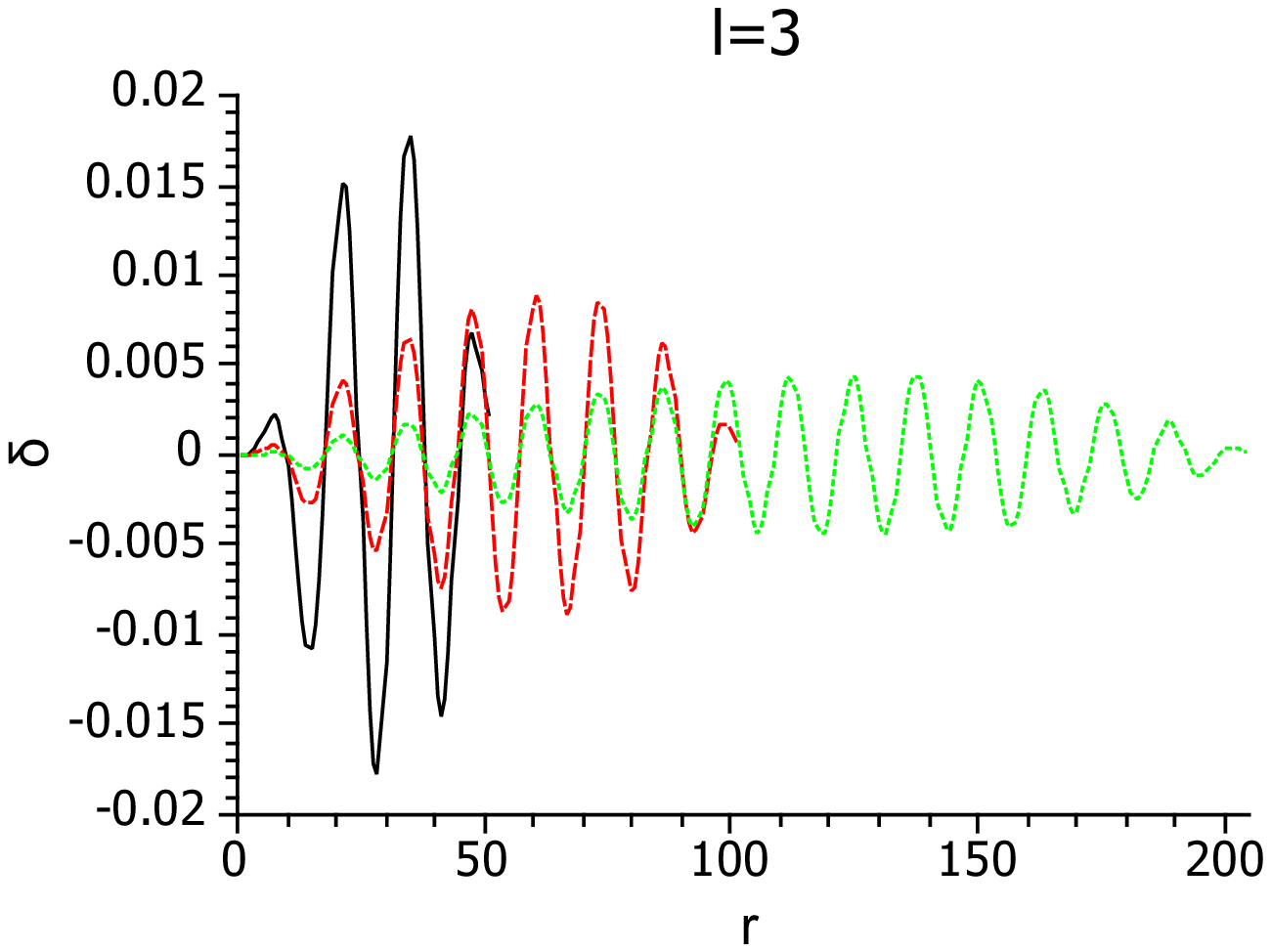}
\includegraphics[angle=-90,width=0.45\columnwidth]{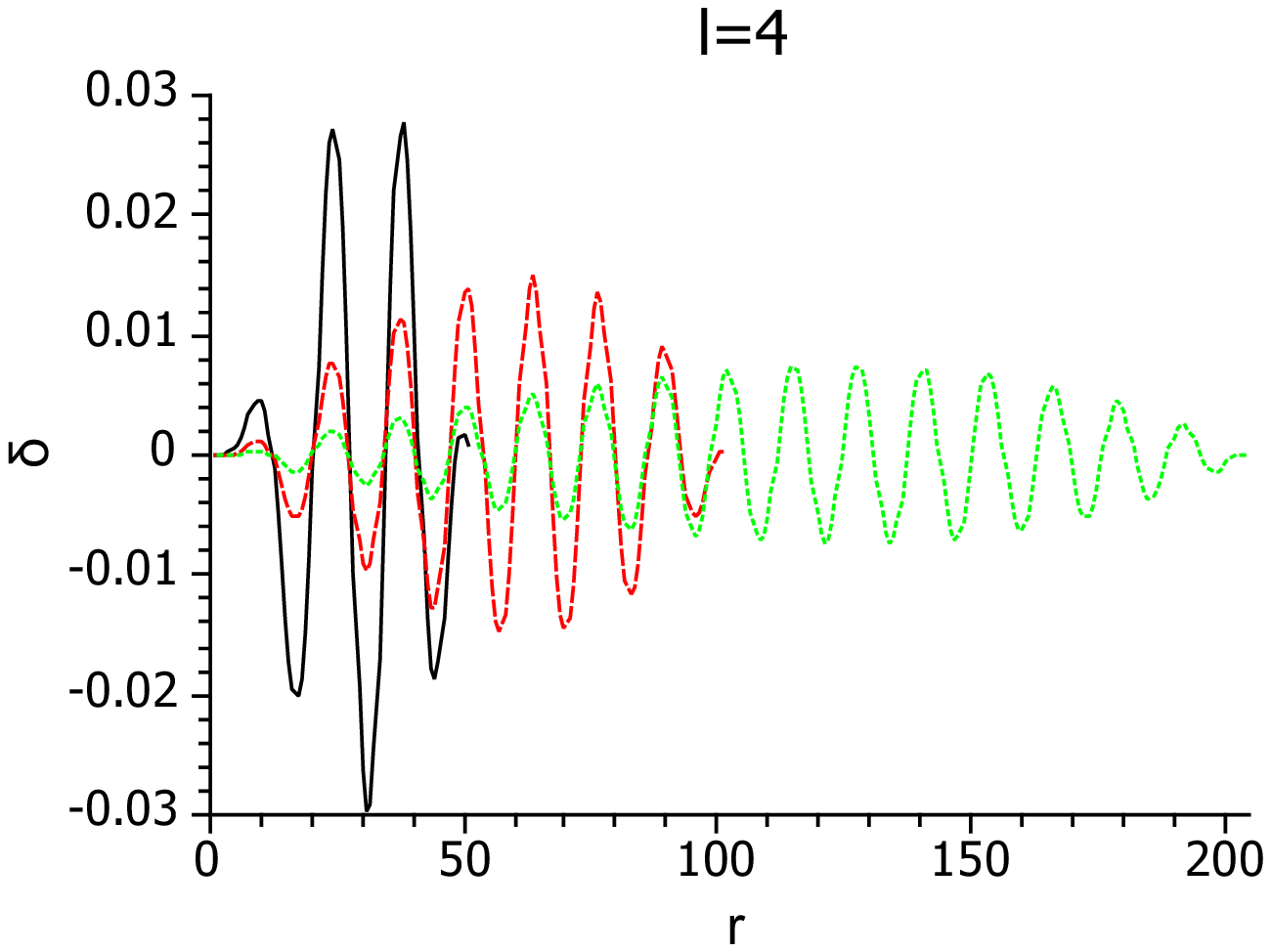}
\caption{Difference between the basis functions and SBFs $\delta=|\tilde{\chi}_{n(k)\ell}(r)-\chi_{\ell}(k_{n\ell}r)|$ for the fixed $k$ = 0.490873843.}
\label{fig:berr_k_corr}
\end{figure}
On Fig.\ref{fig:berr_k_corr}, we plot the $\delta=|\tilde{\chi}_{n(k)\ell}(r)-\chi_{\ell}(k_{n\ell}r)|$, i.e. the difference between the DSBT basis functions and SBFs for the corrected momentum $k_{n\ell}$. It is seen that the difference is an order of magnitude less than in the case of the non-corrected momentum (Fig.\ref{fig:berr}) at the same $r_{max}$. Yet the rate of the basis function convergence to the corrected momentum SBF is still $\Delta k^2$ (whereas $r\ll r_{max}$).
Thus our transformation result arrears to be more exact approximation to the expansion in SBFs on the grid $k_{n\ell}$ than it is on the uniform grid $k_n$. However the grid $k_{n\ell}$ is non-uniform and $\ell$-dependent, and this complicates the application with no convergence rate advantage.

\section{An example: Solution of the 3D time-dependent Schr\"odinger equation} 
\label{Sec:Example}

\subsection{Numerical method}

As an example of DSBT application we developed the method for a numerical solution of 3D time-dependent Schr\"odinger equation (TDSE).


Let us make use of discrete variable representation (DVR). 
We shall begin with the function discretization on the 3D grid in the spherical coordinate system
\begin{eqnarray}
\psi_{ijk}=\Psi(r_i,\arccos\eta_j,\phi_k)r_i\sqrt{\Delta r \Delta\eta_j \Delta\phi}
\end{eqnarray}
Here  $\phi_k=\Delta\phi (k-1),\; k=1,\ldots,N_\phi$ is a polar angle grid; $\Delta\phi=2\pi/N_\phi$ is a grid step; $\eta_j$ and $\Delta\eta_j,\; j=1,\ldots,N_\theta$ are Gauss-Legendre quadrature nodes and weights correspondingly.
Upon DVR implementation, the transformation given by Eqs.(\ref{YProj},\ref{BesselProj}) is written as 
\begin{eqnarray}
\mx{c}=\mx{B}\mx{Y}\bpsi. \label{BYpsi}
\end{eqnarray}
Here $\mx{c}$ is a vector of coefficients of the expansion in spherical waves
\begin{eqnarray}
c_{nlm}=c_{lm}(\tilde{k}_n)\sqrt{w_n}.
\end{eqnarray}
Since further we will need identical momentum grids for all the angular momenta $\ell$, whereas the grid $k_n$ defined by Eq.\eqref{k_grid} differs for odd and even $\ell$, we shall now introduce a new grid for the momentum radial component
\begin{eqnarray} 
\tilde{k}_n = n \Delta k; \quad n=1,\ldots,N, \label{tilde_k_grid}
\end{eqnarray}
and also shall assume $c_{lm}(\tilde{k}_N)\sqrt{w_N}=c_{0lm}$ for even $l$. Validity of this procedure may be justified as follows. Due to the subsection \ref{Sec:Completion}, the coefficient having $n=0$ for all $l$ corresponds to a projection onto a non-regular high-frequency function. In turn, a non-regular high-frequency function may be approximated by the sum of regular functions with large momenta. So, as the highest spectrum part coefficients are evaluated rather inaccurately, one may assume without loss of accuracy the coefficient with $n=0$ to be the value of the projection on a basis function with large momentum $\tilde{k}_N$.
Meanwhile, for a smooth $\Psi(r,\theta,\phi)$ this coefficient vanishes anyway.

The rest of designations used in Eq.\eqref{BYpsi} are as follows:
$\mx{B}$ is the SBT matrix with elements
\begin{eqnarray}
B_{nlmil'm'} &=& [\mx{T}_l\mx{F}_l]_{ni}\delta_{ll'}\delta_{mm'}
\end{eqnarray}
and $\mx{Y}$ is the matrix of the transform to the expansion in spherical harmonics
\begin{eqnarray}
Y_{ilmi'jk}=\delta_{ii'}[\mx{P}\mx{\Phi}]_{lmjk}
\end{eqnarray}
which might be represented as the product of the matrix with elements
\begin{eqnarray}
\Phi_{jmj'k}=\delta_{jj'} \frac{e^{im\phi_k}}{\sqrt{2\pi}}\sqrt{\Delta\phi} 
\end{eqnarray}
and the matrix
\begin{eqnarray}
P_{lmjm'}= \overline{P}_{l}^m(\eta_j)\sqrt{\Delta\eta_j}\, \delta_{mm'}.
\end{eqnarray}
Here $\overline{P}_{l}^m(\eta)$ are associated Legendre polynomials orthonormalized on the Gauss-Legendre quadrature \cite{Melezhik1999}. The transform Eq.\eqref{BYpsi} requires $O(N_r N_\theta N_\phi^2) + O(N_r N_\phi N_\theta^2) + O(N_\theta N_\phi N_r \log_2 N_r )$ operations. It should be noted that the term $O(N_r N_\phi N_\theta^2)$ is caused not only by the operation of polar angle $\mx{P}$ multiplication by the transform matrix, but also by the multiplication by matrices $\mx{T}_l$ (matrix multiplication at the fixed $l$ requires $O(N_rN_\phi l)$ operations, whereas the number of different $l$'s is equal to $N_\theta$). Therefore if one had tried using (instead of DVR) any methods that do not employ $\mx{P}$ transformation, it would not make sense, because it would not imply getting rid of the operations number quadratic growth with $N_\theta$ increasing.  However, the transform algorithm is easily parallelizable, hence may be run in quite modest amount of computer time even for the large value of $N_\theta$.

Besides, we introduce the matrix 
\begin{eqnarray}
\tilde{Y}_{nlmn'jk}=\delta_{nn'}i^l[\mx{P}\mx{\Phi}]_{lmjk}
\end{eqnarray}
of the transform to the expansion in terms of modified spherical harmonics \cite{VMH} related to the common spherical harmonics as $\tilde{Y}_{lm}(\theta,\phi)=i^lY_{lm}(\theta,\phi)$ ($i$ is the imaginary unit here). 
It can be used to perform the switch to the momentum DVR
\begin{eqnarray}
\bphi=\mx{\tilde{Y}}^\dag\mx{B}\mx{Y}\bpsi, \label{YBYpsi}
\end{eqnarray}
where the vector $\bphi$ components relate to the plane wave expansion coefficients as follows 
\begin{eqnarray}
\varphi_{njk}=\varphi(\tilde{k}_n,\arccos\eta_j,\phi_k)\sqrt{\Delta k \Delta\eta_j \Delta\phi}.
\end{eqnarray}
The very possibility of the transition to the plane wave expansion is the main reason for us having introduced the unified grid for the momentum radial components, Eq.\eqref{tilde_k_grid}. 

One may write the Hamiltonian for a particle in external field as 
\footnote{Here and below all equations are expressed in atomic units} 
\begin{eqnarray}
\hat{H}=\frac{\hat{\vec{p}}^{\,2}}{2}-\vec{\mathcal{A}}(t)\hat{\vec{p}}+U(r,\theta,\phi,t).
\end{eqnarray}
Here $\hat{\vec{p}}=-i\nabla$ is the momentum operator, $U(r,\theta,\phi,t)$ is the potential of electron-nuclei interaction, and $\vec{\mathcal{A}}(t)$ is the vector potential of the external electric field. The vector potential definition
\begin{eqnarray}
\vec{\mathcal{A}}(t)=-\int_0^t q\vec{\mathcal{E}}(t')dt,
\end{eqnarray}
slightly differing from the commonly used one, will be used further for the sake of expressions brevity. Here $q$ is a particle charge and $\vec{\mathcal{E}}(t)$ is the external electric field strength. 

Employing DVR and the transform Eq.\eqref{BYpsi} allows to represent the Hamiltonian as a matrix 
\begin{eqnarray}
\mx{H}(t)=\mx{Y}^\dag\mx{B}^\dag[\mx{K}-\mx{\tilde{Y}}(\vec{\mathcal{A}}(t)\vec{\mx{P}})\mx{\tilde{Y}}^\dag]\mx{B}\mx{Y}+\mx{U}(t).\label{H_DVR}
\end{eqnarray}
Here the kinetic energy operator matrix $\mx{K}$, the potential energy operator matrix $\mx{U}$ and the momentum operator matrix $\vec{\mx{P}}$ are diagonal, and their elements are written as:
\begin{eqnarray}
K_{nlmn'l'm'}&=&\frac{k_{nl}^2}{2}\delta_{nn'}\delta_{ll'}\delta_{mm'};\label{KinMat}\\
U_{ijki'j'k'}(t)&=&U(r_i,\theta_j,\phi_k,t)\delta_{ii'}\delta_{jj'}\delta_{kk'};\\
\vec{P}_{njkn'j'k'}&=&\tilde{k}_n\vec{n}_{jk}\delta_{nn'}\delta_{jj'}\delta_{kk'}.
\end{eqnarray}
Hence the operations number needed for the multiplication by matrix \eqref{H_DVR} scales as the one for the transform Eq.\eqref{BYpsi}.
In the absence of the vector potential the multiplication by the Hamiltonian requires only two transforms, namely the direct and inverse ones. When the vector potential is non-zero, one employs the additional couple of angular transforms $\mx{Y}$. In the present case they perform the transition from the expansion in spherical waves to the one in plane waves (that is, to the DVR in the momentum representation) and vice versa. Since all the transforms are orthogonal, the Hamiltonian matrix preserves the original Hamiltonian hermiticity.

Next, we take the opportunity to introduce in the kinetic energy matrix Eq.\eqref{KinMat} the angular momentum-dependent momentum $k_{nl}$. One can define it either just as $k_{nl}=\tilde{k}_n$, or in a more advanced fashion, namely via Eq.\eqref{kcorr}. We will compare these two ways below.

The TDSE has the form
\begin{eqnarray}
i\frac{\partial\psi(r,\theta,\phi,t)}{\partial t}=\hat{H}\psi(r,\theta,\phi,t). \label{TDSE}
\end{eqnarray}
In the matrix form it is written as
\begin{eqnarray}
i\frac{\partial\bpsi(t)}{\partial t}=\mx{H}(t)\bpsi(t). \label{TDSEmatrix}
\end{eqnarray}

Since one may perform the fast fast multiplication by matrix of the form of \eqref{H_DVR}, the time dependent equation can be solved by means of different approaches, e.g. the leap-frog method or short iterative Lanczos propagator method \cite{ParkLight1986}. However here we shall use the split-operator method \cite{Fleck1988}. It might be represented in the form of equations in the following order:
\begin{eqnarray}
\bpsi_1(t)&=&\exp[-i\mx{U}(t+\tau/2)\tau/2]\bpsi(t);\nn\\
\mx{c}_1(t)&=&\mx{B}\mx{Y}\bpsi_1(t);\nn\\
\mx{c}_2(t)&=&\exp[-i\mx{K}\tau/2]\mx{c}_1(t);\nn\\
\bphi_2(t)&=&\mx{\tilde{Y}}^\dag\mx{c}_2(t);\nn\\
\bphi_3(t)&=&\exp[i(\vec{\mathcal{A}}(t+\tau/2)\vec{\mx{P}})\tau]\bphi_2(t);\\
\mx{c}_3(t)&=&\mx{\tilde{Y}}\bphi_3(t);\nn\\
\mx{c}_4(t)&=&\exp[-i\mx{K}\tau/2]\mx{c}_3(t);\nn\\
\bpsi_4(t)&=&\mx{Y}^\dag\mx{B}^\dag\mx{c}_4(t);\nn\\
\bpsi(t+\tau)&=&\exp[-i\mx{U}(t+\tau/2)\tau/2]\bpsi_4(t).\nn
\end{eqnarray}
This sequence of steps is equivalent to $\bpsi(t+\tau)=\exp[-i\mx{H}(t+\tau/2)\tau]\bpsi(t)$ with the accuracy 
$O(\tau^3)$, so the method has a global error of $O(\tau^2)$. As the matrices $\mx{K}$, $\mx{U}$ and $\vec{\mx{P}}$ are diagonal, the exponential functions of them reduce to the exponential functions of complex numbers. Therefore each step of the method performing requires a number of operations $O(N_r N_\theta N_\phi^2) + O(N_r N_\phi N_\theta^2) + O(N_\theta N_\phi N_r \log_2 N_r)$.

Upon the employing of the approach that is being presented, the evolution of the phases of free spherical waves is evaluated more precisely, the greater the evaluation region. It emerges to be an important advantage in comparison to another space approximation techniques that are commonly used today (finite-difference method, finite-element method and so on). This is extremely helpful for the problems that require the consideration of the long-duration wavefunction evolution in large space regions in variable external fields.

It is worth mentioning that, although we are considering the TDSE solving only, the reduction of a problem to the multiplication by the matrix of the form of \eqref{H_DVR} might be also used in iteration methods for the stationary elliptic equations solving as well.

\subsection{Test on 3D time-dependent harmonic oscillator}

\begin{figure}[ht]
\includegraphics[angle=-90,width=0.45\columnwidth]{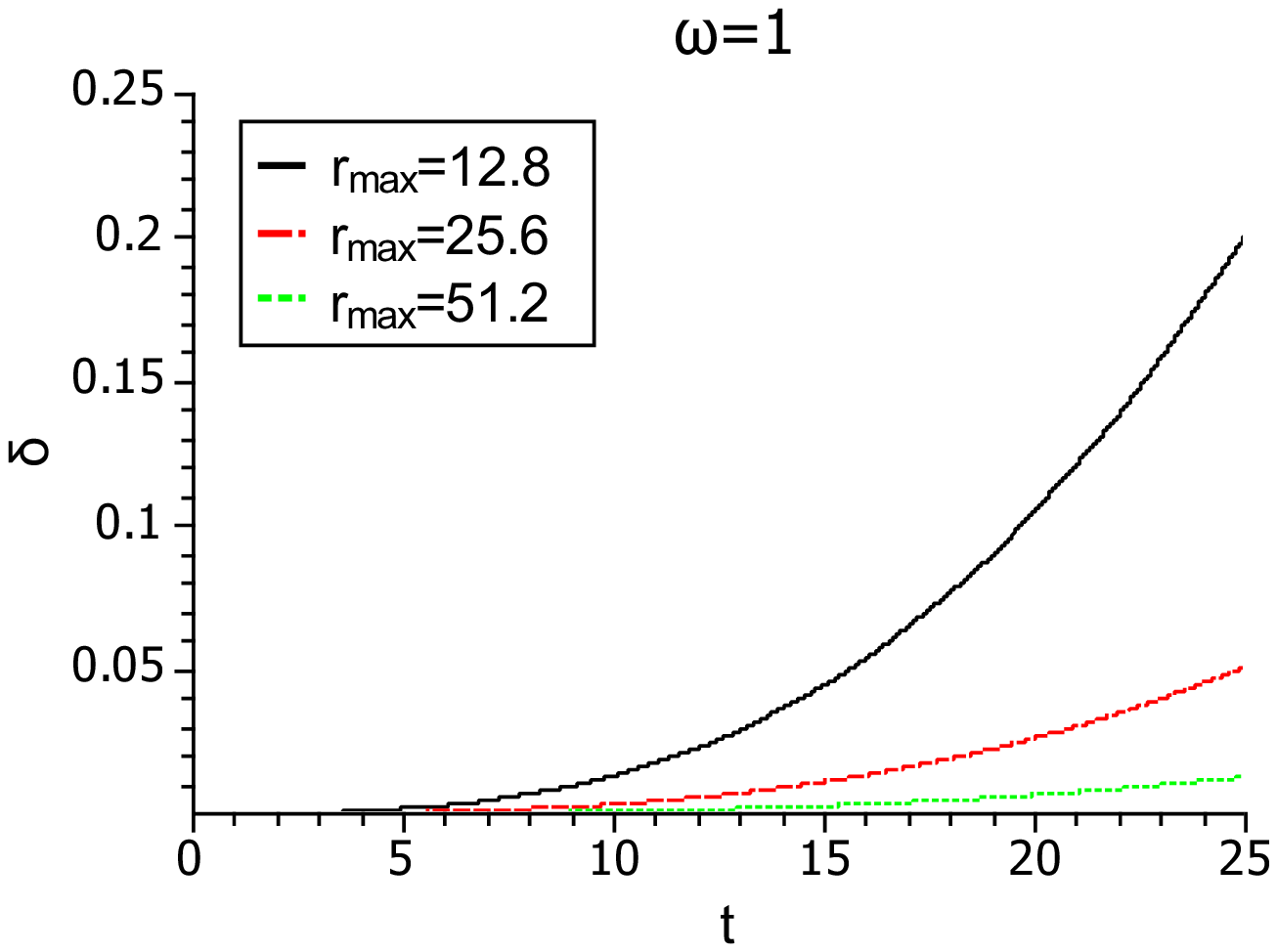}
\includegraphics[angle=-90,width=0.45\columnwidth]{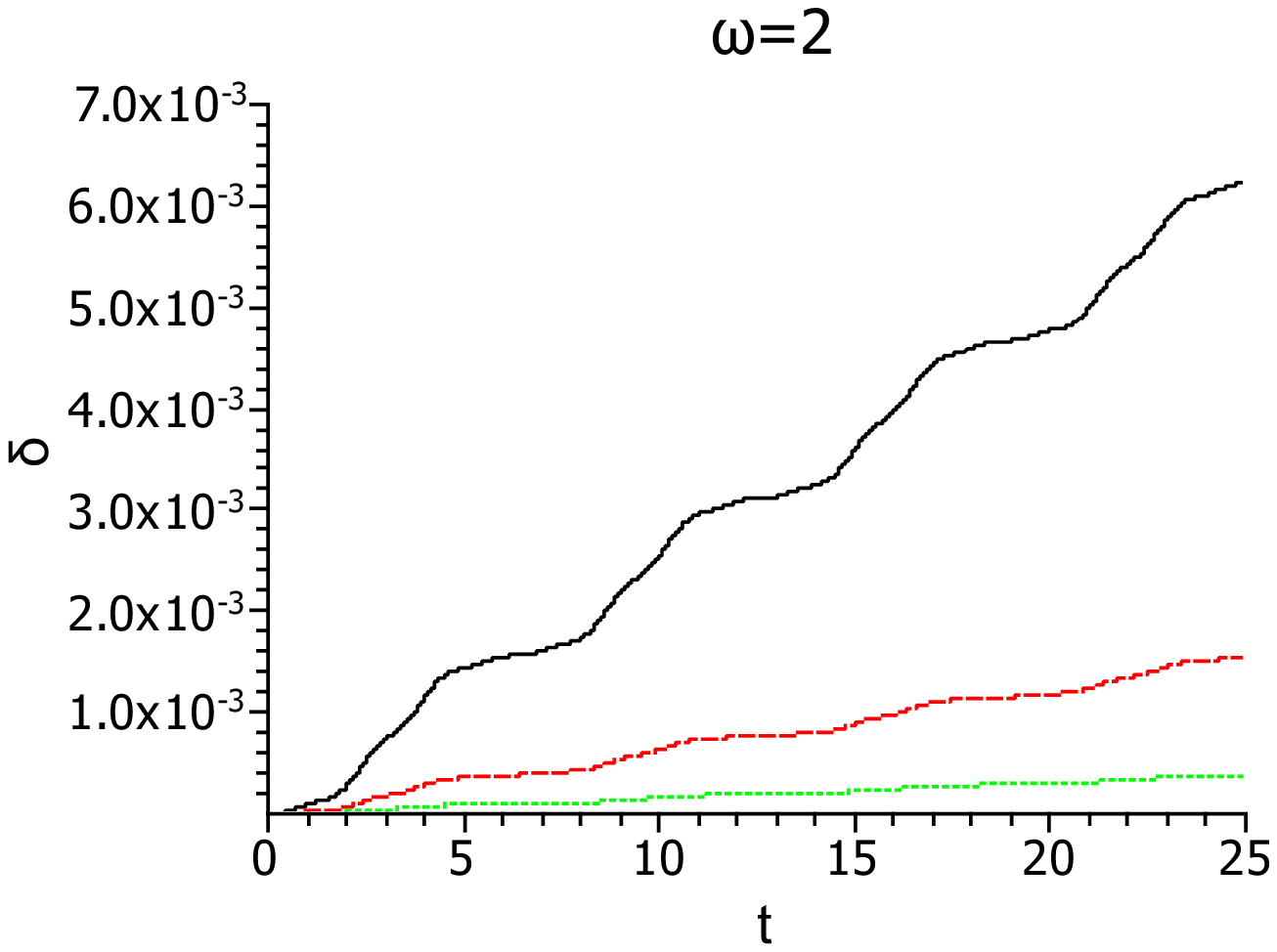}
\\
\includegraphics[angle=-90,width=0.45\columnwidth]{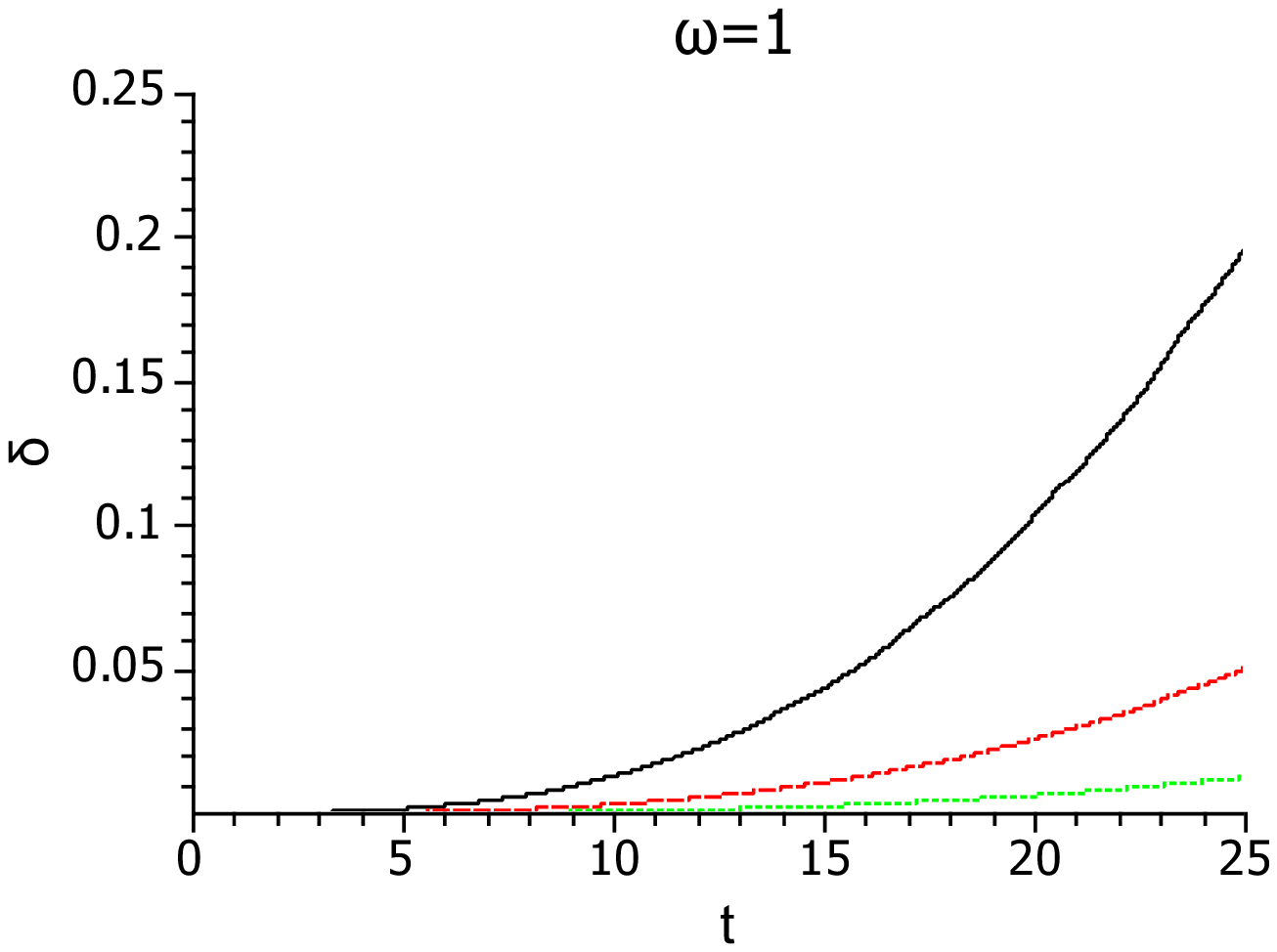}
\includegraphics[angle=-90,width=0.45\columnwidth]{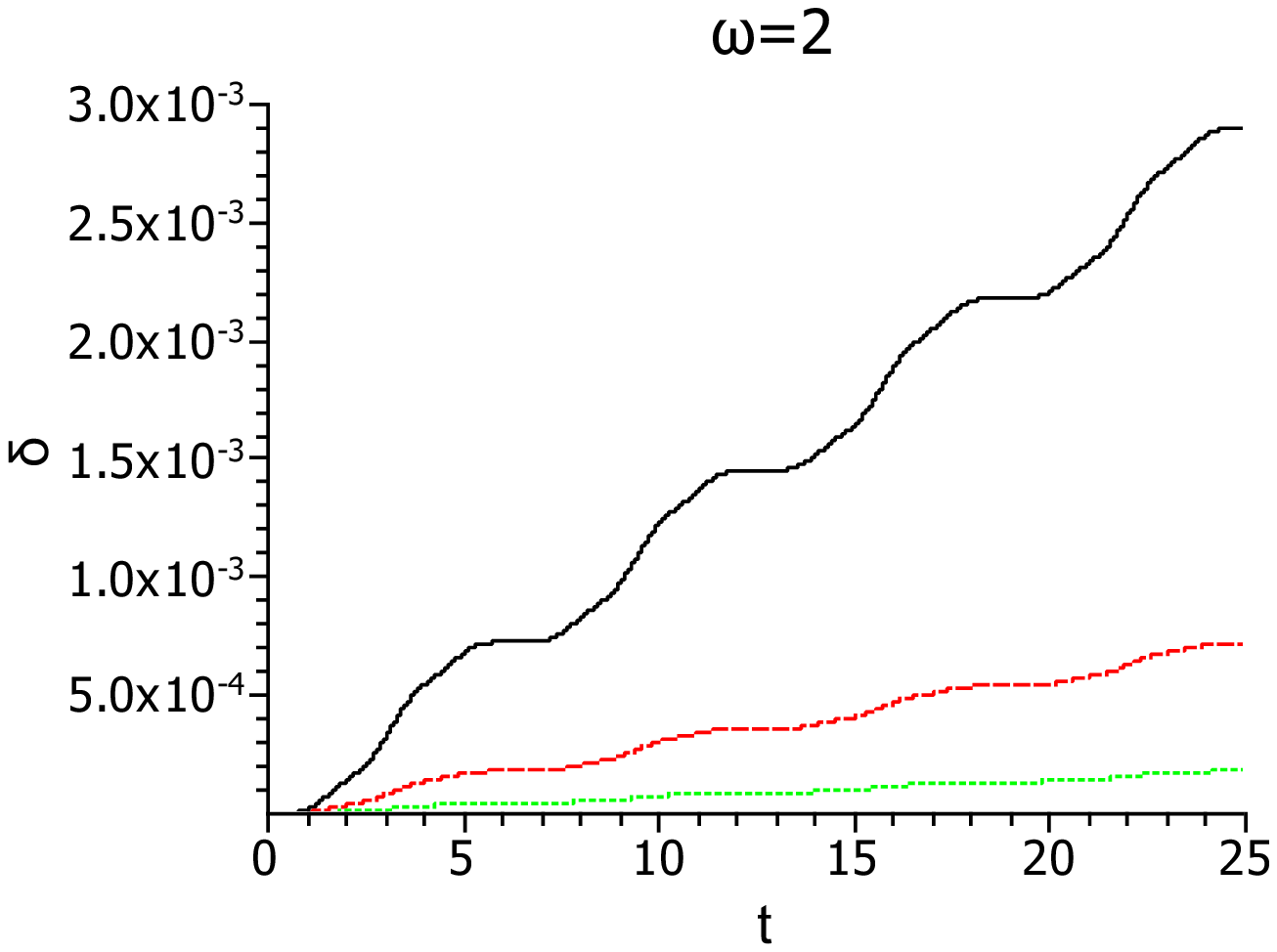}
\caption{Convergence of the solution to the exact one for the oscillator in external time-dependent field: first row --- the coordinate gauge; second row --- the velocity gauge; left column --- external field frequency $\omega=1$; right column --- $\omega=2$.}
\label{fig:Oscill}
\end{figure}
As a first benchmark application let us consider the problem of 3D harmonic oscillator in external field, that has the analytical solution. As this problem possesses features somewhat opposite to those which are optimal for the employing of DSBT-DVR (that is, it needs only the small spatial region size), it proves to be the most stringent test for the method.

The spherically-symmetric three-dimensional harmonic oscillator potential is known to be
\begin{eqnarray*}
U_0(r)=\frac{r^2}{2}.
\end{eqnarray*}
The time-dependent external field can be presented by means of various ways which are equivalent in terms of theory, but different in terms of their implementation by a numerical scheme. We have accomplished the calculations for the external field representation both in the coordinate gauge
\begin{eqnarray*}
&&V(\vec{r},t)=-q\mathcal{E}(t)z;\\
&&\vec{\mathcal{A}}(t)=0.
\end{eqnarray*}
and the velocity gauge
\begin{eqnarray*}
&&V(\vec{r},t)=0;\\
&&\vec{\mathcal{A}}(t)=\mathcal{A}(t)\vec{e}_z.
\end{eqnarray*}
The pulse form was supposed to be 
\begin{eqnarray*}
\mathcal{A}(t)=-\mathcal{A}_0\sin\omega t,
\end{eqnarray*}
and
\begin{eqnarray*}
q\mathcal{E}(t)=-\frac{\partial\mathcal{A}}{\partial t}=\omega\mathcal{A}_0\cos\omega t.
\end{eqnarray*}
We took the external field amplitude to be $\mathcal{A}_0=0.25$. Computations were carried out for the two external field frequencies, namely $\omega=1$ corresponding to the oscillator resonant frequency, and the non-resonant $\omega=2$. 
At the resonant frequency the amplitude of the wavepacket center position oscillations (against center of coordinate)
grows linearly with time, that is, higher and higher spherical harmonics are excited, whereas in the non-resonant case only small $\ell$ harmonics are excited.
Initial state function was set equal to the ground state function. 

We used the parameter $\delta(t)=|1-\langle\psi_\text{osc}(\vec{r},t)|\psi(\vec{r},t)\rangle|$ to estimate the approximation error. Here wave function $\psi_\text{osc}(\vec{r},t)$ is the analytical solution for the three-dimensional harmonic oscillator in a time-dependent external field, and $\psi(\vec{r},t)$ is the numerical solution. 
We set the angular basis parameters $N_\theta=16$ and $N_\phi=1$. In order to diminish the error of the split-operator method (which is of no current interest), the time step has been set very small.

The Figure \ref{fig:Oscill} shows the error $\delta(t)$ of the obtained numerical solution as a function of time at $k_{nl}=\tilde{k}_n$ for the three grids having the same step $\Delta r = 0.2$, but different $\rmax=$ 12.8, 25.6, and 51.2. It is apparent that the solution error falls down as $O(\rmax^{-2})$, as one should expect basing on the fact that the approximate transform error makes the most significant contribution to the overall scheme error at such scheme parameters.

\begin{figure}[ht]
\includegraphics[angle=-90,width=0.45\columnwidth]{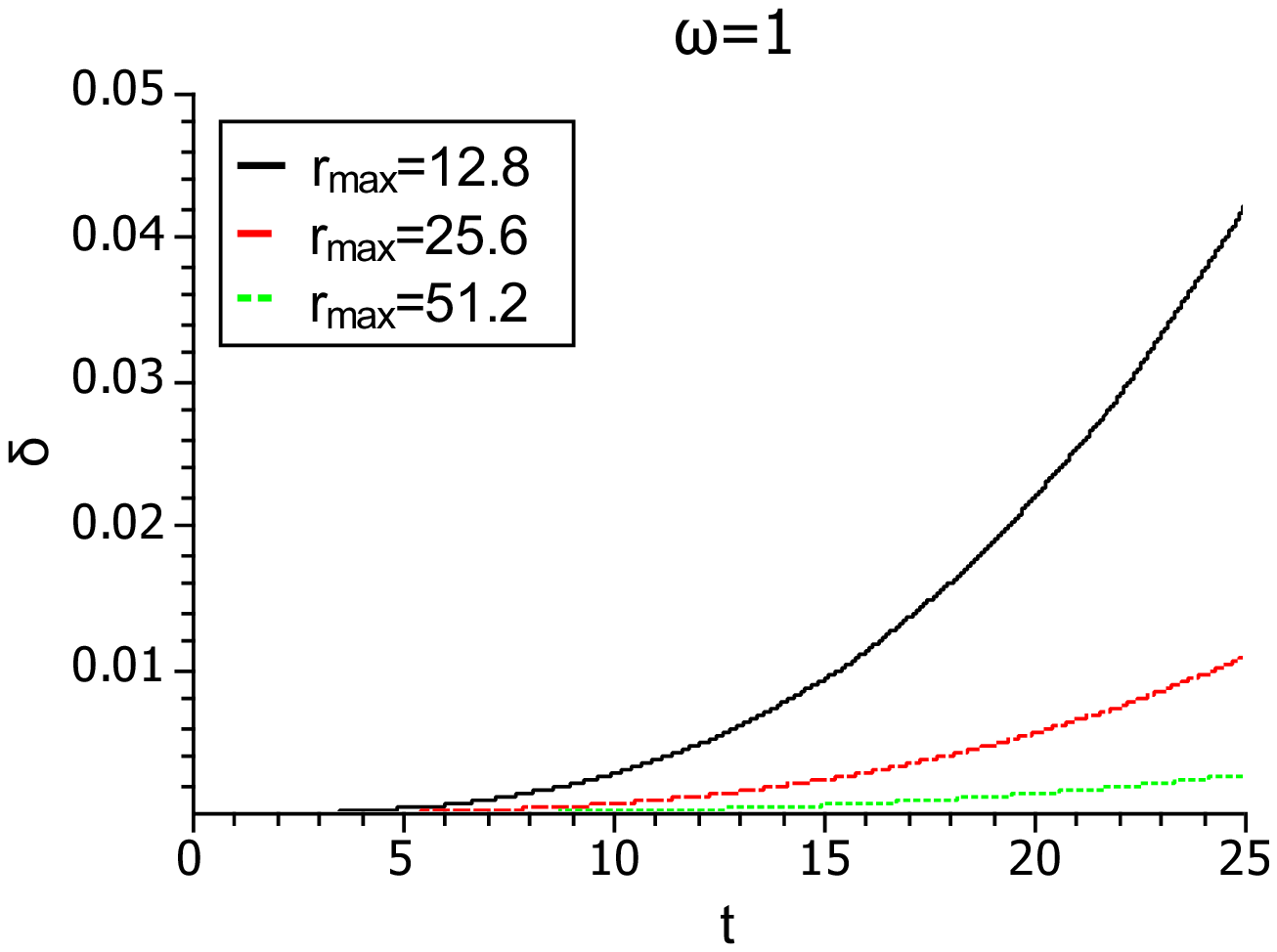}
\includegraphics[angle=-90,width=0.45\columnwidth]{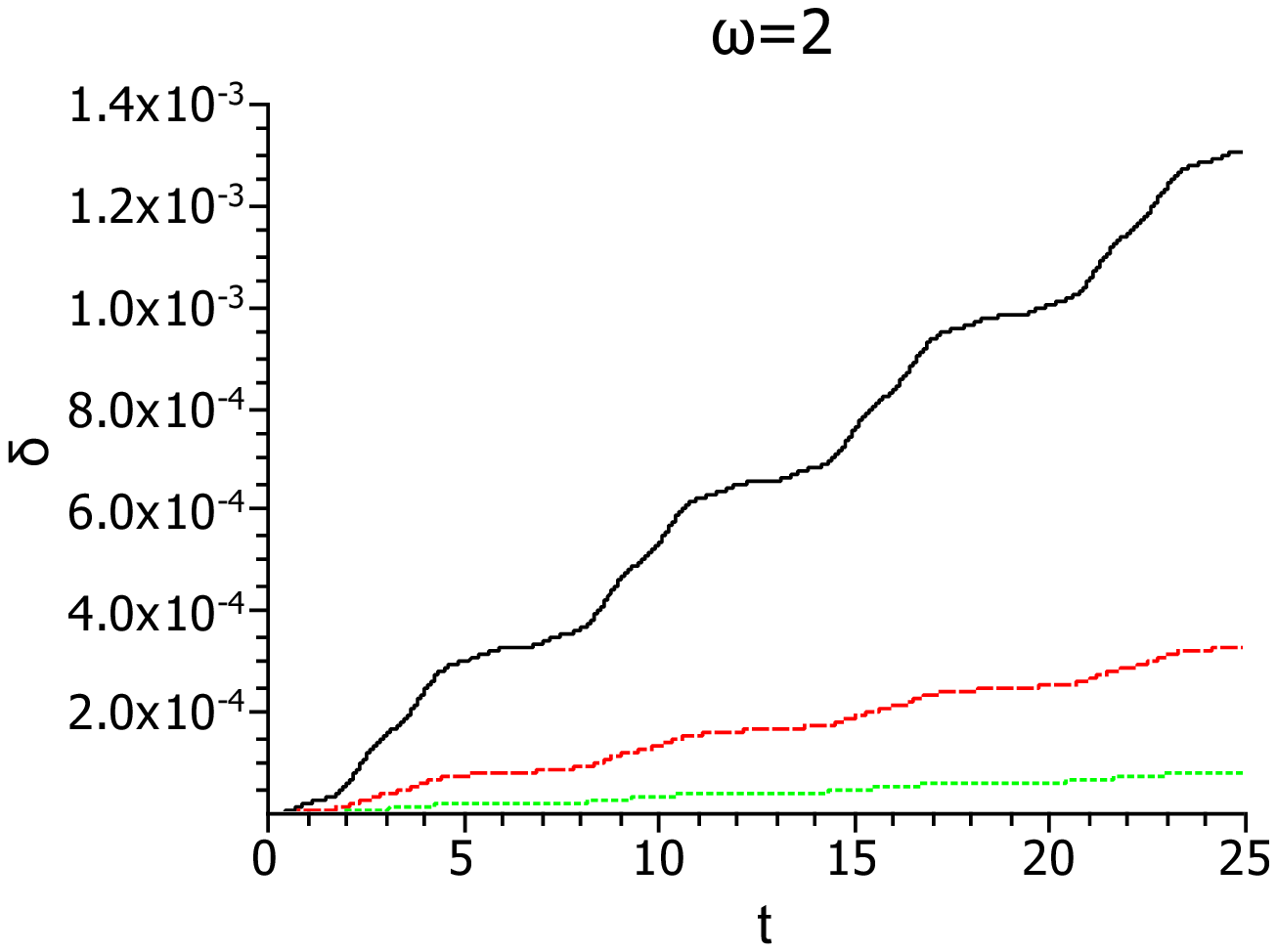}
\\
\includegraphics[angle=-90,width=0.45\columnwidth]{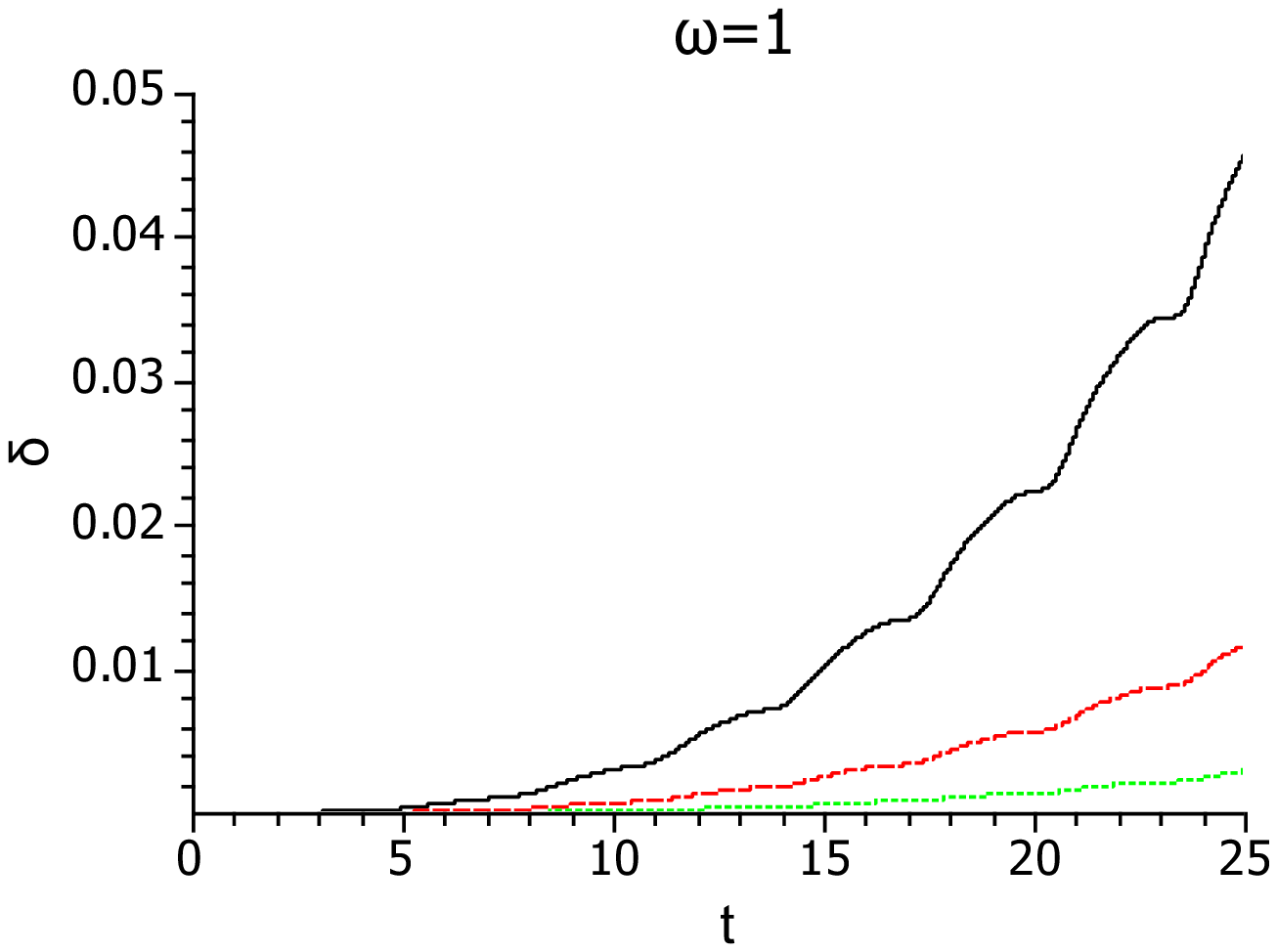}
\includegraphics[angle=-90,width=0.45\columnwidth]{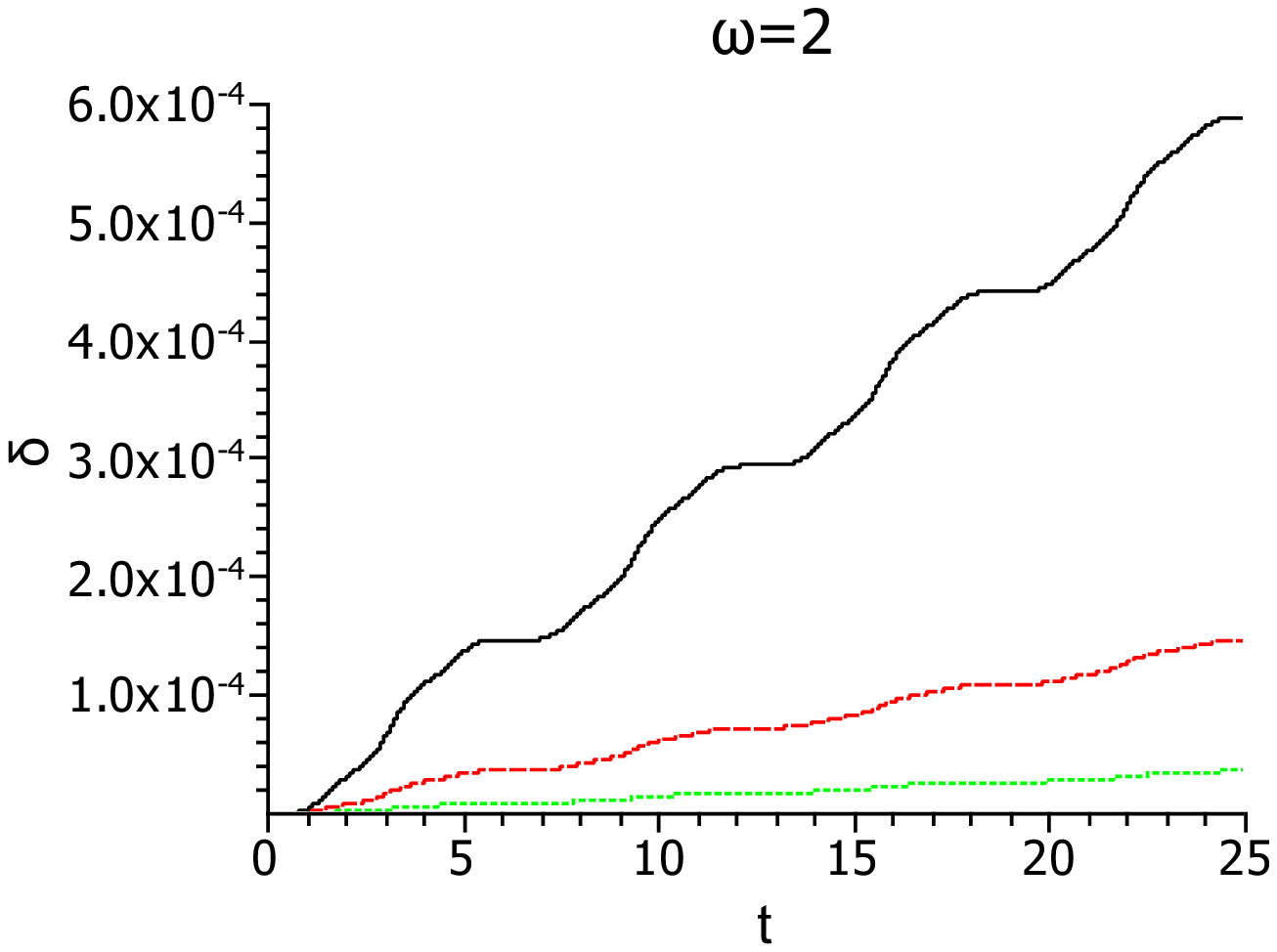}
\caption{Convergence of the solution to the exact one for the oscillator in external time-dependent field for the corrected $k_{nl}$: first row --- the coordinate gauge; second row --- the velocity gauge; left column --- external field frequency $\omega=1$; right column --- $\omega=2$.}
\label{fig:OscillCorr}
\end{figure}
The Figure \ref{fig:OscillCorr} presents the same as the Figure \ref{fig:Oscill}, except that we have employed the corrected $k_{nl}$ from Eq.\eqref{kcorr}. One can see that the rate of the solution convergence to the exact one depending on $\rmax$ is the same as in the case $k_{nl}=\tilde{k}_n$, but the error absolute value is 4 to 5 times less.

\subsection{Molecule in the pump-probe field}

Now let us turn to a problem of more physical use. As a benchmark example we shall consider the H$_2^+$ molecule in the field of complex-shaped laser pulse consisting of the short ultraviolet (XUV) pulse combined with the long infrared radiation (IR) pulse. This emerges as a model of rapidly developing pump-probe techniques \cite{Krausz2014}. The electron is emitted after being subjected to the XUV 
pumping pulse and then moves under joint action of the long-range Coulomb  field and slowly changing IR probe pulse field. Modeling of this process requires computations for a long atomic time period as well as for the large simulation region size $\rmax$ (in order for the electron not to escape outside its boundaries). Since a large $\rmax$ implies a small momentum step, an DSBT based approach emerges to be perfectly appropriate for this problem solving.

First we need to estimate the accuracy that our scheme provides for the singular potential problems which are frequently encountered in atomic physics. To this end, we have computed the eigenenergies of the approximate Hamiltonian \eqref{H_DVR} for the different singular potentials. The ground state energy and wavefunction have been evaluated by means of the imaginary time evolution method (that is  to substitute $t\to -it$ in Eq.\eqref{TDSEmatrix}). The excited states have been evaluated via the imaginary time evolution method with the ortogonalization of the wavefunction to the lower states functions on each time step.

We shall begin with the considering of the Hydrogen atom whose nucleus potential is known to be
\begin{eqnarray*}
U_0(r)=-\frac{1}{r}.
\end{eqnarray*}
\begin{table*}[ht]
\caption{Bound states energies for H.}
\centering
\begin{tabular}{lllll}
\hline\hline
         & & $\Delta r$ & \\
				 \cline{2-4}
 $n$ $\ell$ &	0.2 & 0.1 & 0.05 & Exact \\
\hline
1 0      &	-0.505927 & -0.501575 & -0.500405 & -0.5\\
2 0      &	-0.125738 & -0.125197 & -0.125051 & -0.125 \\
3 0      &	-0.055774 & -0.055614 & -0.055571 & -0.055555(5)\\
2 1      &	-0.125017 & -0.125017 & -0.125017 & -0.125 \\
3 1      &	-0.055572 & -0.055572 & -0.055572 & -0.125 \\
3 2      &	-0.055606 & -0.055606 & -0.055606 & -0.055555(5) \\
\hline\hline
\end{tabular}
\label{tab:HboundsEr102}
\end{table*}
The table \ref{tab:HboundsEr102} demonstrates the convergence of the calculated energy with the grid step $\Delta r$ decreasing at the fixed $r_{max}=102.4$. It is seen that calculated energies of $\ell=0$-states converge to the exact ones quadratically. Meanwhile, $\ell>0$-states energies hardly depend on $\Delta r$ and possess much less errors. The latter are caused mainly by the very SBT error and decrease quadratically with $r_{max}$ increasing. For $\ell=0$, the large error value and rather slow $\Delta r$-convergence result from the fact that the Coulomb wavefunctions with $\ell=0$ have the first derivative discontinuity at $r=0$ and are poorly approximated by the $\sin$ Fourier expansion. 

In order to enhance the convergence rate for $\ell=0$, one can replace the exact Coulomb potential with an effective potential constructed in such a manner that, at a given approximate kinetic energy operator $\mx{Y}^\dag\mx{B}^\dag\mx{K}\mx{B}\mx{Y}$, the approximated Hamiltonian ground state function and energy would coincide with the exact ones for the ground state of a hydrogen-like ion with a nucleus charge $Z$, that is, correspondingly,
\begin{eqnarray}
\varphi_{Z100}(\vec{r})=\frac{Z^{3/2}}{\sqrt{\pi}}\exp(-Zr);\;\; E_{Z10}=-\frac{Z^2}{2}.
\end{eqnarray}
For a nucleus residing in a point with the coordinates $\vec{r}_a$, such a potential is expressed as
\begin{eqnarray}
\tilde{u}_{Z}(\vec{r}_{ijk},\vec{r}_a)=\frac{[\mx{Y}^\dag\mx{B}^\dag[\mx{K}-E_{Z10}\mx{I}]\mx{B}\mx{Y}\bphi_Z(\vec{r}_a)]_{ijk}}{[\bphi_Z(\vec{r}_a)]_{ijk}}, \label{tilde_uZ}
\end{eqnarray}
where
\begin{eqnarray}
[\bphi_Z(\vec{r}_a)]_{ijk}=\varphi_{Z100}(\vec{r}_{ijk}-\vec{r}_a)r_i\sqrt{\Delta r \Delta\eta_j \Delta\phi}.
\end{eqnarray}
However, since $\varphi_{Z100}(\vec{r})$ tends to zero exponentially at large $r$'s, the expression \eqref{tilde_uZ} would yield the result going to infinity at large $|\vec{r}-\vec{r}_a|$ due to numerical errors. To avoid this, we have chosen to use the following potential
\begin{eqnarray}
u_{Z}(\vec{r}_{ijk},\vec{r}_a)=f(|\vec{r}_{ijk}-\vec{r}_a|)\tilde{u}_{Z}(\vec{r}_{ijk},\vec{r}_a)-\frac{[1-f(|\vec{r}_{ijk}-\vec{r}_a|)]Z}{|\vec{r}_{ijk}-\vec{r}_a|}. \label{uZ}
\end{eqnarray}
Here $f(r)$ is the mask function possessing the properties $f(0)=1$, $f(r\to\infty)=0$. In our calculations the mask function of the form
\begin{eqnarray}
f(r)=\exp(-Zr).
\end{eqnarray}
was employed. In such a way, the potential \eqref{uZ} coincides with the potential \eqref{tilde_uZ} when $|\vec{r}-\vec{r}_a|$ is small and with the usual Coulomb potential when it is large. As all the Coulomb functions at $r\to 0$ have an asymptotic behavior  $\sim 1-Zr+O(r^2)$, the increasing of the accuracy of near-$r=0$ approximation of the ground state function $\varphi_{Z100}(\vec{r})$ should lead to the increasing of the accuracy of the approximation of the other Coulomb functions.

\begin{table*}[ht]
\caption{Bound state energies for H with the effective potential in use.}
\centering
\begin{tabular}{llll}
\hline\hline
         & & $\Delta r$ & \\
				 \cline{2-4}
 $n$ $\ell$ &	0.2 & 0.1 & 0.05 \\
\hline
1 0      &	-0.500967 & -0.500133 & -0.500017 \\
2 0      &	-0.125125 & -0.125017 & -0.125002 \\
3 0      &	-0.055593 & -0.055561 & -0.055556 \\
2 1      &	-0.125016 & -0.125017 & -0.125017 \\
3 1      &	-0.055572 & -0.055572 & -0.055572 \\
3 2      &	-0.055606 & -0.055606 & -0.055606 \\
\hline\hline
\end{tabular}
\label{tab:corrHboundsEr102}
\end{table*}
The table \ref{tab:corrHboundsEr102} exhibits the same as the Table \ref{tab:HboundsEr102} does, except that the potential
\begin{eqnarray*}
U_0(r)=u_{Z}(\vec{r},0)
\end{eqnarray*}
has been used instead of the Coulomb one. One can easily observe the decreasing of the differences between the calculated energies of $\ell=0$-states and the exact energies of H atom stationary states, whereas for the $\ell>0$-states these differences apparently do not increase.

Now turn to the molecular Hydrogen ion H$_2^+$. When dealing with multinuclear systems, one has to employ a potential with singularities that do not coincide with the coordinate origin. We shall write the approximate nuclear potential in the hydrogen molecule in the following way:
\begin{eqnarray*}
U_0(r)=u_{1}(\vec{r},\vec{R}/2)+u_{1}(\vec{r},-\vec{R}/2).
\end{eqnarray*}
where $\vec{R}$ is the internuclear vector with the length $R=2$ which corresponds to the equilibrium internuclear distance for the H$_2^+$ ground state. We have chosen the internuclear direction along the $Oz$ axis orientation, $\vec{R}||\vec{e}_Z$.
\begin{table*}[ht]
\caption{Bound state energies for H$_2^+$.}
\centering
\begin{tabular}{lllll}
\hline\hline
       &    & $N_\theta$ & \\
			 \cline{2-4}
 State &	4 & 8 & 16 & Exact \\
\hline
 $1\sigma_g$ &	-1.066449 & -1.094991 & -1.101242 & -1.102634 \\
 $2\sigma_u$ &	-0.618383 & -0.659020 & -0.666117 & -0.667534 \\
\hline\hline
\end{tabular}
\label{tab:H2plusboundsEdr02}
\end{table*}
The table \ref{tab:HboundsEr102} manifests the calculated energies for H$_2^+$ ground and first excited states converge with the angular basis parameter $N_\theta$ increasing at the fixed $\Delta r = 0.2$ and $r_{max}=102.4$. The ``exact'' energies given here were obtained through the calculation via the method \cite{Serov2002} based upon the spheroidal coordinates utilizing. The error arising from the grid step is negligible in this case, therefore the table of convergence over the grid step is not presented here.

Next let us consider the evolution of the molecular ion H$_2^+$ in the field of two overlapping linearly polarized laser pulses
\begin{eqnarray*}
\vec{\mathcal{A}(t)}= \mathcal{A}_{UV}(t)\vec{n}_{UV} +  \mathcal{A}_{IR}(t)\vec{n}_{IR}.
\end{eqnarray*}
Here the ``XUV'' pulse was supposed to have the Gaussian envelope
\begin{eqnarray*}
\mathcal{A}_{UV}(t)=-\mathcal{A}_{UV}\exp\left(-2\ln 2 \frac{t^2}{w_{UV}^2}\right)\cos\omega t
\end{eqnarray*}
where $w_{UV}$ is the full width at half maximum .
Next, the ``IR'' pulse was chosen to have a compact support and the $\cos^2$--envelope, as follows
\begin{eqnarray*}
\mathcal{A}_{IR}(t)=-\mathcal{A}_{IR} \cos^2[\pi (t-t_{IR})/\tau_{IR}] \cos\omega_{IR}(t-t_{IR}),\, |t-t_{IR}|<\tau_{IR}/2,
\end{eqnarray*}
where $\tau_{IR}$ is the overall pulse duration, $t_{IR}$ is the shift of the arrival time of the IR-pulse center relative to that for the XUV pulse. The external field of this form is employed in the attosecond streaking method \cite{Krausz2014}. The XUV-pulse triggers the ionization, then the detected electrons spectrum dependence on the time shift $t_{IR}$ enables to determine the IR pulse genuine form, or, in the case of this form being known, to obtain the time delay of the electron emission during the ionization process.

The probe pulse parameters was taken to be $\omega_{IR}=0.062832$, $\mathcal{A}_{IR}=0.05$, and $\tau_{IR}=2T_{IR}=200$ (which are common values in modern attosecond streaking experiments), and the pump pulse parameters, correspondingly, were $\omega_{UV}=|E_0|+0.5$ ($E_0$ standing for the molecule ground state energy, evaluated by means of the imaginary time evolution method) $\mathcal{A}_{UV}=0.25$, and $w_{UV}=10$. Both pulses polarization were chosen to be co-directed with the molecular axis, $\vec{n}_{UV}=\vec{n}_{IR}=\vec{e}_Z$. In all the examples referred to below we used the numerical scheme parameters as follows: $r_{max}=409.6$, time step $\tau=\Delta r^2/4$, evolution beginning time $t_0=-\tau_{IR}/2+t_{IR}$, and evolution termination time $t_{fin}=\tau_{IR}/2=100$.

\begin{figure}[ht]
\includegraphics[angle=-90,width=0.45\columnwidth]{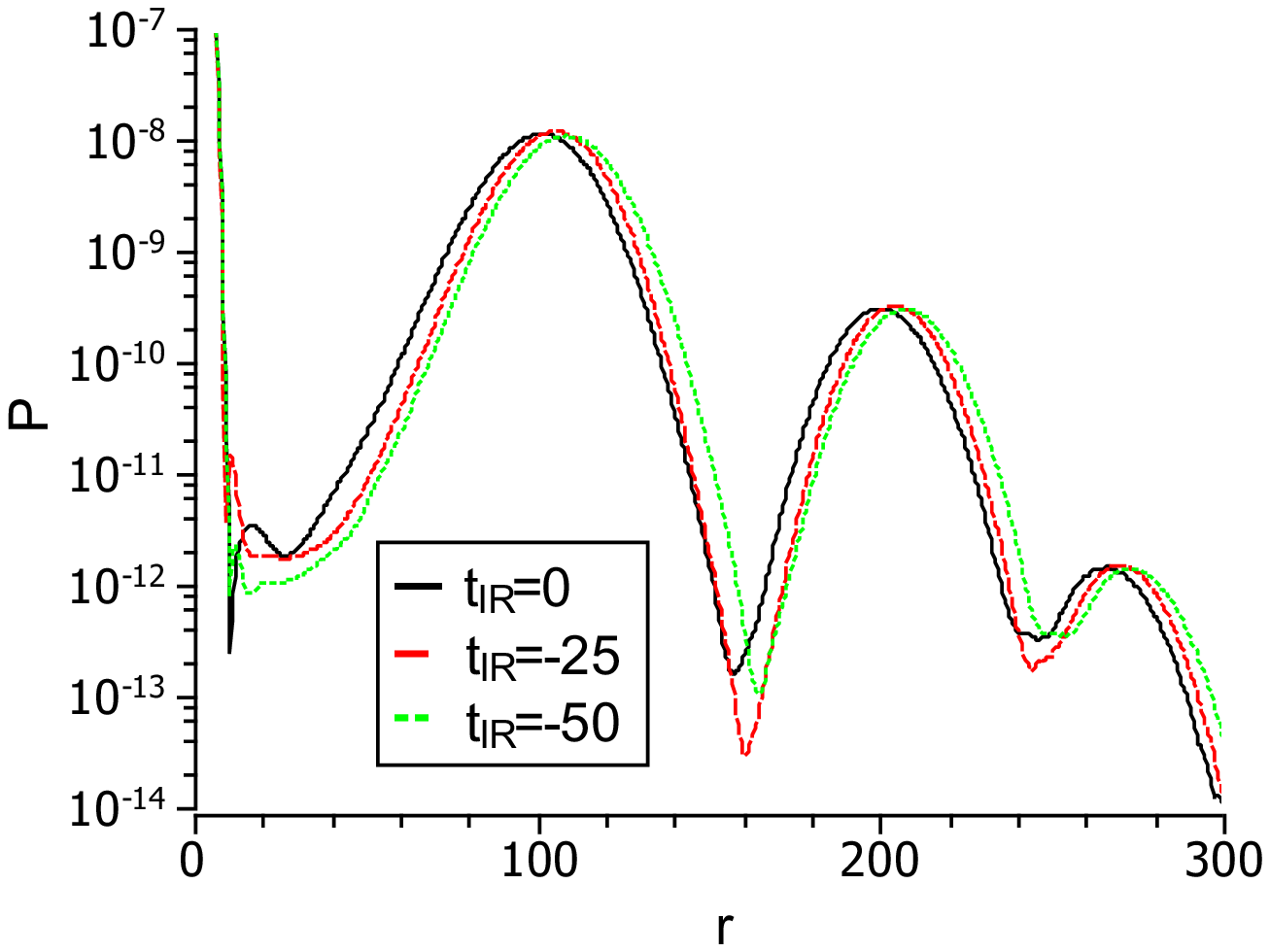}
\includegraphics[angle=-90,width=0.45\columnwidth]{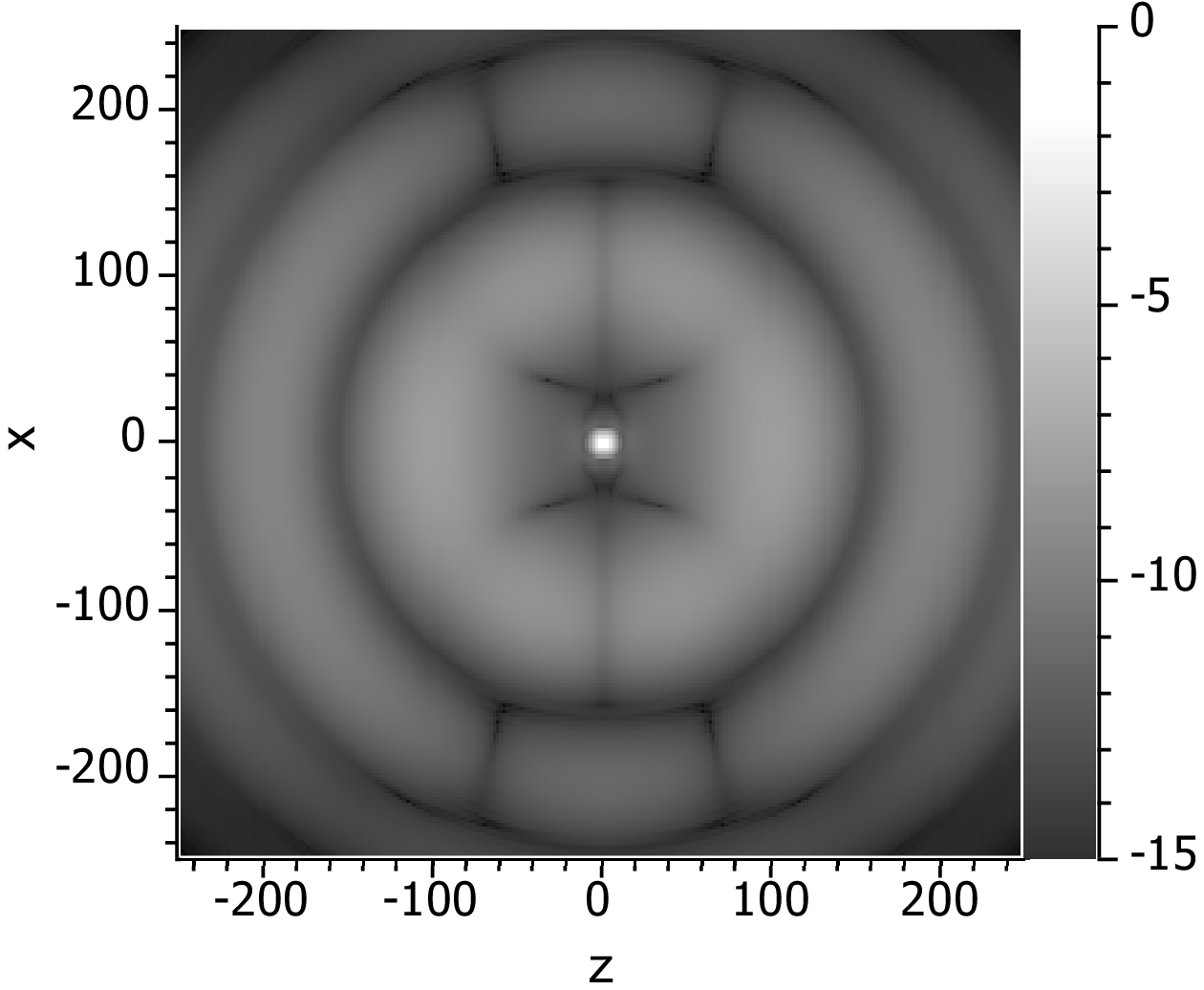}
\caption{(left) The probability density $P=|\Psi(\vec{r},t_{fin})|^2$ versus $r$ at $\theta=0$ for the different $t_{IR}$'s; (right) gray scale map of $\log_{10} P$ versus Cartesian coordinates $x$ and $z$ at $y=0$ for $t_{IR}=0$.}
\label{fig:H2plusAS}
\end{figure}
The Fig. \ref{fig:H2plusAS} shows the probability density $P(\vec{r},t_{fin})=|\Psi(\vec{r},t_{fin})|^2$ that the electron is at $\vec{r}$. 
The calculations were performed for the scheme parameters $N_\theta=16$ and $\Delta r = 0.2$. Due to the stationary phase approximation, for the time $t_{fin}\gg w_{UV}$ and for $r\gg 1$, the relation $P(\vec{r},t_{fin})\sim \sigma(\vec{r}/t_{fin})$ holds, where $\sigma(\vec{k})$ is the differential cross section of electrons emission depending on momentum. On the left panel of the Fig. \ref{fig:H2plusAS}, the peak near $r=0$ corresponds to the wavefunction of the ground and other stationary states, whereas the peak centered in the vicinity of $r=100$ emerges due to the one-photon ionization, and and the rest large $r$ peaks are caused by the multiphoton processes. This is apparently confirmed by the right panel of the Fig. \ref{fig:H2plusAS}, where the $r=100$ enhanced probability ring
has one node depending on $\theta$, the circle of larger radius has  two nodes corresponding to the dipole and quadrupole distributions arising from the absorption of one or two photons correspondingly. The left panel of the Figure \ref{fig:H2plusAS} also demonstrates the probability density dependence on the IR pulse phase at the moment of the XUV pulse arrival. The theory predicts the probe pulse action causing the electron momentum shift equal to the magnitude of the IR pulse at the moment of the electron emission from the molecule (which roughly coincides with the moment of the XUV pulse arrival). This is exactly what is observed on the right panel of the Figure \ref{fig:H2plusAS}.

\begin{figure}[ht]
\includegraphics[angle=-90,width=0.45\columnwidth]{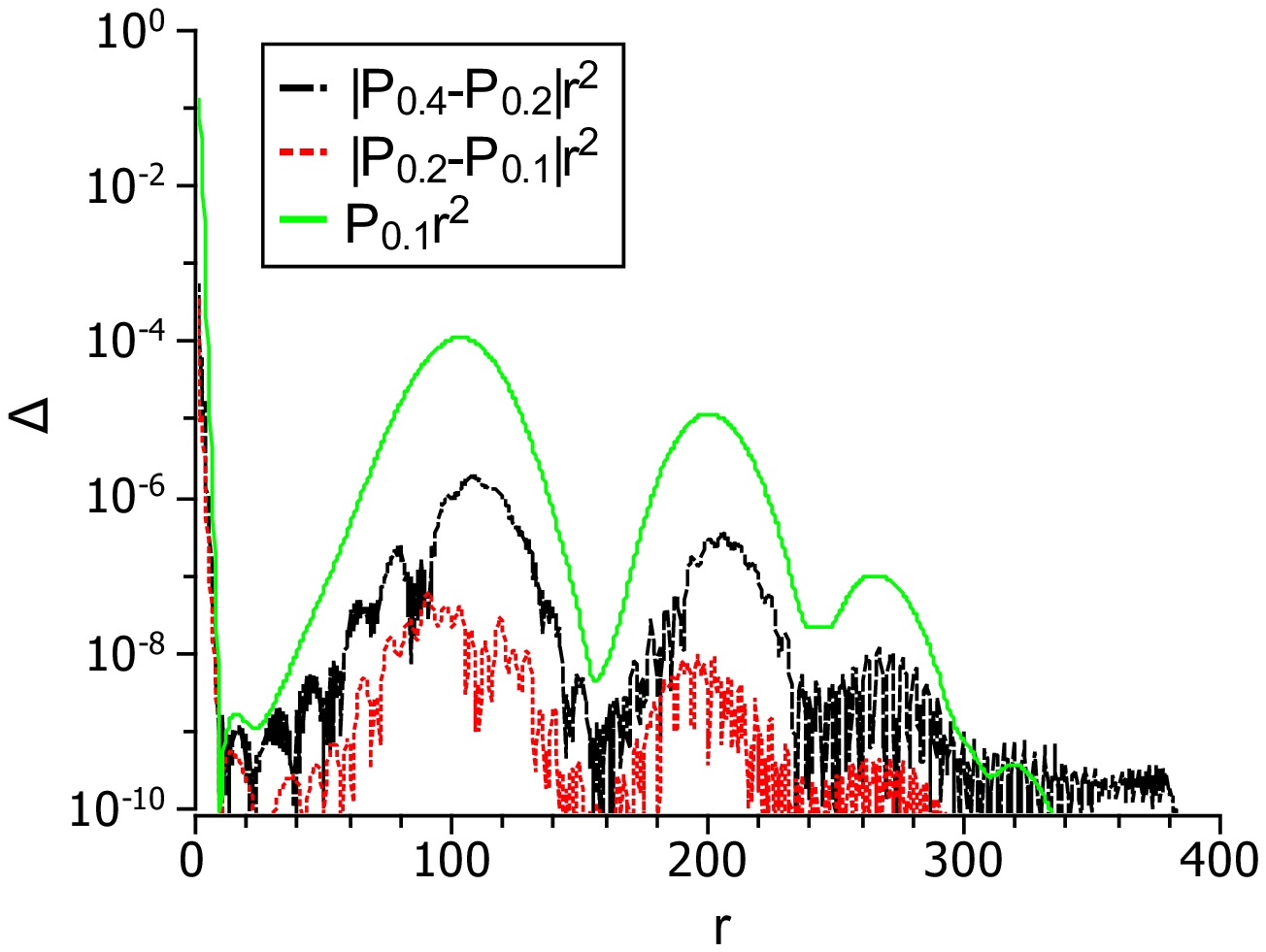}
\includegraphics[angle=-90,width=0.45\columnwidth]{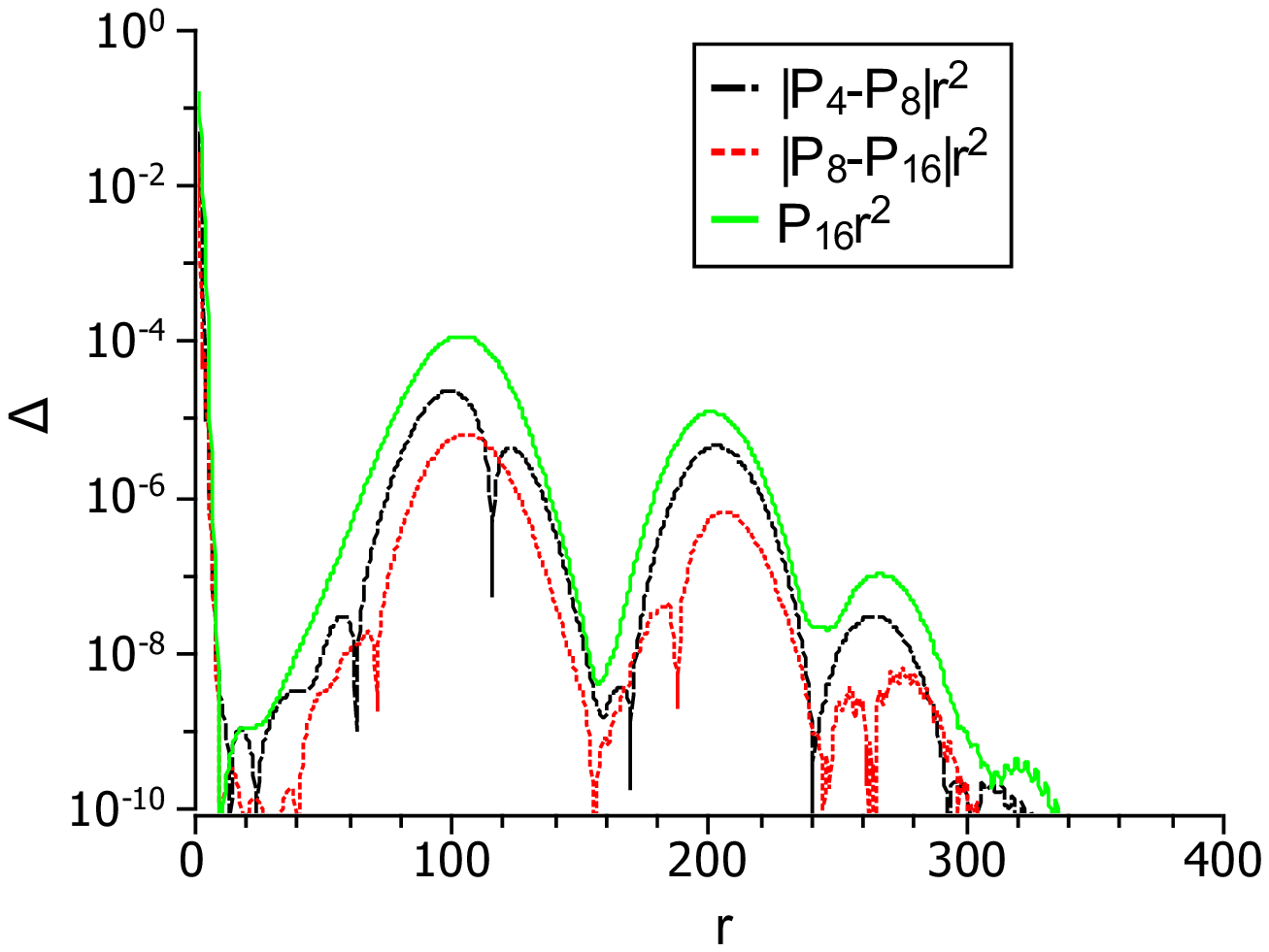}
\caption{(left) The pairwise differences $\Delta=|P_{2\Delta r,N_\theta}-P_{\Delta r,N_\theta}|r^2$ at the fixed $N_\theta=8$; (right) The pairwise differences $\Delta=|P_{\Delta r,N_\theta/2}-P_{\Delta r,N_\theta}|r^2$ at the fixed $\Delta r=0.2$.}
\label{fig:H2plusASconv}
\end{figure}
Besides, we have examined the probability density $P$ convergence rate depending on the step $\Delta r$ and on the angular basis size, $N_\theta$. The left panel of the \ref{fig:H2plusASconv} presents the pairwise differences $\Delta(r)=|P_{2\Delta r,N_\theta}(r\vec{e}_Z,t_{fin})-P_{\Delta r,N_\theta}(r\vec{e}_Z,t_{fin})|r^2$ for the probability densities $P_{\Delta r,N_\theta}(\vec{r},t_{fin})$ evaluated on the three grids having the steps $\Delta r = $ 0.4, 0.2 and 0.1 at the fixed $N_\theta=8$. For the sake of comparison $P_{\Delta r,N_\theta}(r\vec{e}_Z,t_{fin})r^2$ for $\Delta r = $ 0.1 and $N_\theta=8$ is plotted on the same figure. It is apparent that even on the coarsest grid with $\Delta r = $ 0.4 the error is of the order of 1\%; this value is quite small in terms of experimental accuracy which is common in the field in question.
Upon halving the step size, the error drops down by 1-2 orders of magnitude. However, the error in a particular point decreases non-uniformly, actually as expected due to the global basis functions using.

The right panel of the \ref{fig:H2plusASconv} displays the pairwise differences $\Delta(r)=|P_{\Delta r,N_\theta/2}(r\vec{e}_Z,t_{fin})-P_{\Delta r,N_\theta}(r\vec{e}_Z,t_{fin})|r^2$ for the three different angular bases with $N_\theta=$ 4, 8, and 16 at the fixed step $\Delta r = $ 0.2, as well as $P_{\Delta r,N_\theta}(r\vec{e}_Z,t_{fin})r^2$ for $\Delta r = $ 0.2 and $N_\theta=16$. One can see that the error due to the angular basis small size is much larger than that due to the radius step. This is related to the molecular potential non-centrality. For $N_\theta=4$ the error has magnitude about 25\% (in the vicinity of maxima), whereas upon the basis size increasing up to $N_\theta=8$ the error drops down to 6\%. Therefore $N_\theta=16$ has been chosen for the main part of our calculations.

\section{Conclusion} \label{Sec:Conclusion}

We have developed the algorithm for the DSBT that possesses the advantages of orthogonality, performing fastness and uniform grid. Our approach is based upon the SBT factorization into the two subsequent orthogonal transforms, namely the fast Fourier transform (requiring the operations number $O(N\log_2 N)$) and the orthogonal transform founded on the discrete orthogonal Legendre polynomials (requiring the operations number $O(\ell N)$). Our discrete transform converges to the exact SBT as the square of the momentum grid step. 

Besides, basing on DSBT and DVR, we have also elaborated the 3D TDSE solving method (DSBT-DVR).
The examination of the DSBT-DVR algorithm has demonstrated its efficiency for the purposes of solving of time-dependent problems in atomic and molecular physics. An DSBT based approach allows to evaluate the free spherical wave functions evolution the more accurately, the more is the spatial region size. 
It appears to be an advantage in comparison to another methods applied in this field.
This is especially helpful for problems like the modelling of the attosecond streaking approach and other pump--probe techniques, since they require the computation of the wavefunction evolution under the joint action of long-lasting pulses and the weak Coulomb field on large spatial regions. Another important preference of the method proposed is the fast convergence over grid step when applied to the problems with smooth (or artificially smoothed) potentials.

It should be noted that the current DSBT-DVR version does not make any use of another helpful DBBT feature, namely the DSBT capability to be employed for the aim of the evaluation of multi-center integrals \cite{Toyoda2009}. The leveraging of this capability for the solving of both SSE and TDSE for the multielectron molecules is expected to be the matter of our future work.

\section*{Acknowledgements}
The author thanks Dr. Tatiana Sergeeva for help in the preparing of the text of this paper. Also, the author wish to thank Dr. Serguei Patchkowskii for helpful discussions. The author acknowledges support of the work from the Russian Foundation for Basic Research (Grant No. 14-01-00520-a).

\appendix

\section{Proof of transform orthogonality}

Let us begin with the demonstration of the Eq.\eqref{sumDDPp} validity.
Consider the sum
\begin{eqnarray} 
\sigma_s \equiv \sum_{i=0}^{N} \nabla [P_\ell](i,N-1) i^s
\end{eqnarray}
It may be transformed as follows
\begin{eqnarray*} 
 && \sigma_s = \sum_{i=0}^{N} P_\ell(i,N-1) i^s - \sum_{i=0}^{N} P_\ell(i-1,N-1) i^s \\
 && = \sum_{i=0}^{N-1} P_\ell(i,N-1) i^s + P_\ell(N,N-1) N^s \\
 && - \sum_{i=1}^{N} P_\ell(i-1,N-1) i^s - P_\ell(-1,N-1) 0^s
\end{eqnarray*}
According to \eqref{intDPx}, if $s<\ell$ the sums in the last string equal zero, hence
\begin{eqnarray} 
\sigma(s) = P_\ell(-1,N-1) \left[ (-1)^\ell N^s -  0^s \right]. \label{sigmas}
\end{eqnarray}
Here we also used the DLOP parity property \cite{Neuman1974}
\begin{eqnarray} 
P_\ell(i,N)=(-1)^\ell P_\ell(N-i,N). \label{DLOPparity}
\end{eqnarray}
Now let us take an arbitrary discrete polynomial 
\begin{eqnarray}
p(i)=\sum_{s=0}^\mu C_s i^s; \quad \mu<\ell,
\end{eqnarray}
and consider the sum
\begin{eqnarray} 
\sigma \equiv \sum_{i=0}^{N} \nabla [P_\ell](i,N-1) p(i) = \sum_{s=0}^\mu C_s \sigma_s
\end{eqnarray}
By using Eq.\eqref{sigmas} we obtain
\begin{eqnarray} 
\sigma = P_\ell(-1,N-1) \left[ (-1)^\ell p(N) - p(0) \right] \label{sum1DDPp}
\end{eqnarray}
Next, one can construct the weighted sum
\begin{eqnarray*} 
&& \sum_{i=0}^{N} \nabla [P_\ell](i,N-1) p(i) w_i(N) = \\
&& \sigma -\frac{1}{2}\nabla [P_\ell](0,N-1) p(0) -\frac{1}{2}\nabla [P_\ell](N,N-1) p(N) 
\end{eqnarray*}
Making use of the Eqs.(\ref{sum1DDPp}, \ref{backDiff}, \ref{DLOPparity}) and the normalization condition $P_\ell(0,N-1)=1$ yields
\begin{eqnarray} 
 && \sum_{i=0}^{N} \nabla [P_\ell](i,N-1) p(i) w_i(N) = \nonumber \\
 && \frac{1+P_\ell(-1,N-1)}{2} \left[ (-1)^\ell p(N) - p(0) \right].
\end{eqnarray}
After the division of both sides of this equation by $[1+P_\ell(-1,N-1)]/2$, we arrive to Eq.\eqref{sumDDPp}.

Now let us prove Eq.\eqref{orthoImL}. As this equation is symmetric with respect to the exchange of indices $n$ and $m$, for definiteness we shall assume $n>m$. Since $L_{ml}=0$ for $l>m$, one can write
\begin{eqnarray}
 && \sum_{l=p_\ell}^{N_\ell} [\delta_{nl} - L_{nl}] [\delta_{ml} - L_{ml}] = - L_{nm} + \sum_{l=p_\ell}^{m} L_{nl} L_{ml} \label{orthoImLsec}
\end{eqnarray}
For sake of the notation simplicity, from now on we designate
\begin{eqnarray}
 \lambda_{nm} \equiv \sum_{l=p_\ell}^{m} L_{nl} L_{ml} 
\end{eqnarray}
So, according to Eq.\eqref{orthoImLsec}, Eq.\eqref{orthoImL} holds true, when
\begin{eqnarray}
 \lambda_{nm} = L_{nm}. \label{LLL}
\end{eqnarray}
The substitution of $\mx{L}$ elements definition from Eq.\eqref{Ldef} yields
\begin{eqnarray}
 \lambda_{nm} = \sum_{l=p_\ell}^{m} P_\ell'(n-l,2n) P_\ell'(m-l,2m) \left(1-\frac{\delta_{l0}}{2}-\frac{\delta_{lm}}{2}\right). \label{lambdanm}
\end{eqnarray}
By change of the summation index $i=m-l$ we can rewrite \eqref{lambdanm} as
\begin{eqnarray}
 \lambda_{nm} = \sum_{i=0}^{m-p_\ell} P_\ell'(n-m+i,2n) P_\ell'(i,2m) \left(1-\frac{\delta_{im}}{2}-\frac{\delta_{i0}}{2}\right).  \label{ilambdanm}
\end{eqnarray}
Due to Eq.\eqref{DLOPparity}, DDLOP have the parity property 
\begin{eqnarray} 
 P_\ell'(i,N)=(-1)^{\ell-1} P_\ell'(N-i,N). \label{DDLOPparity}
\end{eqnarray}
The sum in Eq.\eqref{ilambdanm} might be split into the two sums as
\begin{eqnarray}
 \lambda_{nm} &=& \frac{1}{2}\sum_{i=0}^{m} P_\ell'(n-m+i,2n) P_\ell'(i,2m) \left(1-\frac{\delta_{i0}}{2}\right) \nonumber \\
 &+& \frac{1}{2}\sum_{i=0}^{m-1} P_\ell'(n-m+i,2n) P_\ell'(i,2m) \left(1-\frac{\delta_{i0}}{2}\right). \label{s2lambdanm}
\end{eqnarray}
Here we made use of the fact that $P_\ell'(m,2m)=0$ when $p_\ell=1$ at even $\ell$. 
Next we apply the parity property to both DDLOPs in the summand of the second sum in Eq.\eqref{s2lambdanm} and make the summation index change $i\to 2m-i$ 
\begin{eqnarray*}
 && \sum_{i=0}^{m-1} P_\ell'(n-m+i,2n) P_\ell'(i,2m) \left(1-\frac{\delta_{i0}}{2}\right)= \\
 && \sum_{i=0}^{m-1} P_\ell'(n+m-i,2n) P_\ell'(2m-i,2m) \left(1-\frac{\delta_{i0}}{2}\right)= \\
 && \sum_{i=m+1}^{2m} P_\ell'(n+m-i,2n) P_\ell'(i,2m) \left(1-\frac{\delta_{i,2m}}{2}\right)
\end{eqnarray*}
The latter sum then might be combined with the (remained unchanged) first sum in Eq.\eqref{s2lambdanm} to get
\begin{eqnarray}
 \lambda_{nm} = \frac{1}{2}\sum_{i=0}^{2m} P_\ell'(n-m+i,2n) P_\ell'(i,2m) w_i(2m).
\end{eqnarray}
As $P_\ell'(n-m+i,2n)$ is the polynomial of the order $\ell-1$, we can apply Eq.\eqref{sumDDPp} to obtain  
\begin{eqnarray*}
 && \sum_{i=0}^{2m} P_\ell'(n-m+i,2n) P_\ell'(i,2m) w_i(2m) \\ 
 && = (-1)^{\ell} P_\ell'(n+m,2n) - P_\ell'(n-m,2n).
\end{eqnarray*}
Finally, after using the parity property Eq.\eqref{DDLOPparity}, the result becomes
\begin{eqnarray}
 \lambda_{nm} = - P_\ell'(n-m,2n).
\end{eqnarray}
Since $L_{nm}= - P_\ell'(n-m,2n)$ for $n>m\geq n_{0\ell}>0$, we have thus proved Eq.\eqref{LLL} and therefore Eq.\eqref{orthoImL}.


\end{document}